\documentclass[onefignum,onetabnum]{siamart220329}

\usepackage{lipsum}
\usepackage{amsfonts}
\usepackage{graphicx}
\usepackage{epstopdf}
\usepackage{algorithmic}
\usepackage{enumitem}
\usepackage{booktabs}
\usepackage{comment}

\ifpdf
  \DeclareGraphicsExtensions{.eps,.pdf,.png,.jpg}
\else
  \DeclareGraphicsExtensions{.eps}
\fi

\definecolor{red2}{RGB}{204,0,0}
\definecolor{blue2}{RGB}{0,103,165}
\hypersetup{colorlinks,linkcolor=[RGB]{0,103,165},citecolor=[RGB]{180,0,0}}

\newcommand{\RV}[1]{{\color{black}{#1}}}

\newsiamremark{remark}{Remark}
\newsiamremark{hypothesis}{Hypothesis}
\crefname{hypothesis}{Hypothesis}{Hypotheses}
\newsiamthm{claim}{Claim}

\usepackage{dsfont}

\def\be{\begin{equation}}
\def\ee{\end{equation}}

\def\x{\mathbf{x}}

\def\y{\mathbf{y}}

\def\f{\mathbf{f}}
\def\z{\mathbf{z}}
\def\v{\mathbf{v}}

\def\b{\mathbf{b}}

\def\K{\mathbf{K}}

\def\prox{\texttt{prox}}

\newcommand{\norm}[1]{\left\lVert#1\right\rVert}

\headers{Positive polynomial approximations in high dimensions}{Y. Chen, D. Xiu And X. Zhang}

\title{On enforcing non-negativity in polynomial approximations in high dimensions}

\author{Yuan Chen\footnotemark[1]~, Dongbin Xiu\thanks{E-mail addresses: \texttt{\{chen.11050, xiu.16\}@osu.edu}. Department of Mathematics, The Ohio State University, Columbus, OH 43210, USA. Funding: This work was partially supported by AFOSR FA9550-22-1-0011.} \and Xiangxiong Zhang\footnotemark[2]\thanks{E-mail address: \texttt{zhan1966@purdue.edu}. Department of Mathematics, Purdue University, 150 N. University Street, West Lafayette, Indiana 47907. X.Z. was supported by NSF DMS-2208515.}}

\usepackage{amsopn}

\ifpdf
\hypersetup{
  pdftitle={Title},
  pdfauthor={Y. Chen, D. Xiu And X. Zhang}
}
\fi

\usepackage{soul}

\begin{document}

\maketitle

\begin{abstract}
Polynomial approximations of functions are widely used in scientific computing. 
In certain applications, it is often desired to require the polynomial approximation to be non-negative (resp. non-positive), or bounded within a given range, due to constraints posed by the underlying physical problems. Efficient numerical methods are thus needed to enforce such conditions. In this paper, we discuss effective numerical algorithms for polynomial approximation under non-negativity constraints. We first formulate the constrained optimization problem, its primal and dual forms, and then discuss efficient first-order convex optimization methods, with a particular focus on high dimensional problems.
Numerical examples are provided, for up to $200$  dimensions, to demonstrate the
effectiveness and scalability of the methods.
\end{abstract}



\section{Introduction}

We are interested in finding a polynomial approximation to an unknown multivariate function $f(\x)$, via its samples $f(\x_i)$, $i=1,\dots,K$, $\x_i \in \mathbb{R}^d$, $d\geq 1$. In many applications, the function is expected to obey a set of given constraints, usually in the form of inequalities. For examples, density or mass is supposed to be positive, or non-negative. Consequently, it is highly desirable to enforce the polynomial approximations to satisfy the same constraints. Violation of the constraints results in non-physical numerical results and often catastrophic break down of the underlying numerical model.

For polynomial approximations with constraints, most of the existing efforts can be classified into two kinds of approaches: construction-based approach and optimization-based approach. The first approach resorts to characterization of the sought-after polynomials with given structures, such as positivity. Examples of such methods include the positivity-preserving interpolation \cite{butt1993preserving,hussain2008positivity,campos2019algorithms,de2017positive, despres2020computation}, with an extension to bound-preserving \cite{campos2020projection}. These methods usually require precise and complicated algebraic derivations. They usually work with relatively simple constraints and are difficult to apply in high dimensions. The second approach is to enforce the constraints by solving a constrained optimization problem, where the constraints are usually replaced by their approximations in the polynomial space.  Convex optimization tools are leveraged to obtain the solution. \RV{In \cite{allen2022bounds},  Bernstein basis that forms a nonnegative partition of unity was used, where the positivity constraint on the whole polynomial becomes the positivity on the expansion coefficients. The resulting approximation is constrained uniformly. However, this approach is not capable of reproducing polynomials \cite{zala2020structure}.} In \cite{zala2020structure}, the feasible set formed by the constraints is characterized by a set of hyperplanes. A greedy algorithm is then employed to iteratively enforce the constraints represented by the hyperplanes.

Although we do not consider solving partial differential equations (PDEs) in this paper, it is worth mentioning the approaches for positivity-preserving numerical methods for PDEs. It is often crucial for numerical schemes for solving PDEs to satisfy certain constraints for the sake of stability, such as positivity of density and pressure in fluid dynamics equations. Bound-preserving and positivity-preserving numerical schemes include  \cite{zhang2010maximum,zhang2010positivity} for conservation laws, \cite{fan2022positivity,zhang2017positivity} for compressible Navier–Stokes equations, \cite{wu2018positivity} for magnetohydrodynamics, diffusion equations \cite{shen2020unconditionally,ju2021maximum}, \RV{transport equations \cite{peterson2024optimization,bochev2020optimization} and \cite{kirby2024high,barrenechea2024nodally} for more general PDEs}, to name a few.

The focus of this paper is on the development of a numerical framework for finding constrained polynomial approximation of functions, which should be scalable for high dimensions. For simplicity we focus on the non-negative constraint, which can be extended to more general inequality constraints. We formulate the non-negative polynomial approximation problem as a convex minimization with the constraints enforced over a set of finite number chosen points. 
In high dimensions, the cardinality of the polynomial space can be exceedingly large, resulting in  a large scale minimizaation problem.
To circumvent the challenge, we employ first-order convex optimization methods, which use only gradient of the cost functions and thus scales well with the problem size. Moreover, the corresponding dual problem is formulated in a space whose dimension is the number of the constraint points. This is a user's choice and controllable in practice.

 There is a vast amount of literature on first-order optimization methods. In this paper, we investigate a few widely known methods, including Fast Iterative Shrinkage-Thresholding Algorithm (FISTA) \cite{beck2009fast,nesterov2013gradient}, accelerated Primal-Dual Hybrid Gradient (PDHG) method \cite{chambolle2011first,pock2011diagonal}, Douglas–Rachford splitting \cite{lions1979splitting}, which is equivalent to Alternating Direction Method of Multipliers (ADMM) \cite{gabay83} and split Bregman method \cite{Goldstein:2009:SBM:1658384.1658386} with special parameters. 
Moreover, the FISTA algorithm can be accelerated \cite{beck2017first} by using an adaptive restart technique \cite{nesterov2013gradient, o2015adaptive}. 

In this paper, we demonstrate that the FISTA with restart method applied to a proper dual problem is efficient to enforce non-negativity at finite given locations for a polynomial approximation. For simpler problems, it is possible to design nearly optimal parameters to accelerate convergence for methods like the Douglas-Rachford splitting \cite{liu2023simple}. 
For simplicity, we do not consider tuning parameters. We design a comprehensive set of numerical examples, from one dimension to hundreds of dimensions, to compare these popular first-order methods. 
 The numerical cost is estimated and verified in experiments to demonstrate the scalability and efficiency in practical use, which indicates that the restarted FISTA on the dual problem is a simple and effective method without tuning parameters.

The rest of the paper is organized as follows. In Section \ref{sec:setup}, we introduce the setup of the problem. In Section \ref{sec:method}, we provide a detailed description of our method. A comprehensive numerical study is presented in Section \ref{sec:numex}. 
\section{Problem Setup}
\label{sec:setup}
\subsection{Least squares polynomial approximation}

Consider the problem of approximating an unknown function $f: D \mapsto \mathbb{R}$, $D \subset \mathbb{R}^d$, $d \geq 1$. Let $\mathbf{x}=(x_1,...,x_d)$ be the variable and $f\in L^2(D)$. Consider  the  subspace of polynomials of degree up to $n \geq 1$
\be
    \Pi_n^d:=\operatorname{span}\left\{\mathbf{x}^{\mathbf{k}}=x_1^{k_1} \cdots x_d^{k_d},|\mathbf{k}| \leq n\right\},
\ee
where $\mathbf{k}=(k_1,\dots,k_d)$ is multi-index with $|\mathbf{k}|=k_1+\cdots+k_d$. The dimension of this subspace is 
\be
    N=\operatorname{dim} \Pi_n^d=\left(\begin{array}{c}n+d \\d \end{array}\right)=\frac{(n+d) !}{n ! d !}.
\ee
Let $\{\psi_j(\x),j=1,...,N\}$ be an orthonormal basis of $\Pi_n^d$, then any $\tilde f \in \Pi_n^d$  can be expressed as
\be
    \tilde{f}(\mathbf{x})=\sum_{|\mathbf{k}|=0}^n c_{\mathbf{k}} \psi_{\mathbf{k}}(\mathbf{x})=\sum_{k=1}^N c_k \psi_k(\mathbf{x}),
\ee
where we have used a linear ordering to map the multi-index $\mathbf{k}$ to a single index $k$. By using vector notation
\be
    \mathbf{\Psi}(\mathbf{x})=\left[\psi_1(\mathbf{x}), \ldots, \psi_N(\mathbf{x})\right]^{{T}},
\ee 
the expression can be written as 
\be
    \widetilde f(\x) = \langle \mathbf{c},  \mathbf{\Psi}(\x) \rangle,
\ee
where $\mathbf{c} = [c_1,...,c_N]^T$ and $\langle \cdot,\cdot \rangle$ is the dot product.

We then consider approximating the function $f$ using its samples. Let $\x_1,\dots,\x_K,$ be a sequence of sample locations in the domain $D$, $\mathbf{f}: = \left[f(\x_1),...,f(\x_K)\right]^T$ be the sample function values, and
\be
   \mathbf{\Psi}_a =  \begin{bmatrix}  \mathbf{\Psi}^T(\x_1)  \\   \cdots  \\ \mathbf{\Psi}^T(\x_K)  \end{bmatrix} \in \mathbb{R}^{K\times N},
\ee
be the Vandermonde-like matrix, the standard least squares approximation problem can be solved via
\be
    \label{equ:probls1}
    \min_{\mathbf{c}} \norm{\mathbf{\Psi}_a \mathbf{c}-\mathbf{f}}_2^2,
\ee
under the assumption that the problem is over-determined with $K>N$.

\subsection{Constrained Polynomial Approximation}
\label{sec:setup-constraint}

Next we consider a polynomial approximation under linear inequality constraints on the coefficients $\mathbf{c}$:
\be
    \label{equ:prob1}
    \min_{\mathbf{c}} \norm{\mathbf{\Psi}_a \mathbf{c}-\mathbf{f}}_2^2,\qquad \text{subject to}\qquad \mathbf{B} \mathbf{c} \geq \b,
\ee
where $\mathbf{B} \in \mathbb{R}^{C \times N}$, $\b \in \mathbb{R}^C$, $C$ is the number of the constraints, and the inequality is enforced component-wise. 
For a set $\Omega$, the 
  indicator function is defined as
\be \label{indicator}
\iota_{\Omega}(x)= 
\begin{cases}
0, & x \in \Omega, \\ 
+\infty, & x \notin \Omega.
\end{cases}
\ee 
We can rewrite problem \eqref{equ:prob1} in the following equivalent form with a constant $\alpha>0$:
\be 
   \label{equ:primalmin}
   \min_{\mathbf{c} \in \mathbb{R}^N} ~~\frac{1}{2}\alpha\norm{\mathbf{\Psi}_a \mathbf{c}-\mathbf{f}}_2^2 + \iota_{\Lambda}( \mathbf{B}\mathbf{c} - \b), 
\ee
where 
$$\Lambda = \{ \mathbf{a}\in\mathbb{R}^C: \mathbf{a}\geq \mathbf{0}\}.$$
The indicator function term is non-differentiable but convex. The matrix $\mathbf{\Psi}_a^T\mathbf{\Psi}_a$ is positive semidefinite thus the $l^2$ approximation term is also convex.

Enforcing non-negativity at some points can be recasted as the constraint  $\mathbf{B}\mathbf{c} \geq \b$. 
Let $\y_1$,...,$\y_M \in D$ be a sequence of points of interest, and let $\mathbf{\Psi}_p \in \mathbb{R}^{M \times N}$ be the Vandermonde-like matrix at these points,
    \be
         \mathbf{\Psi}_p =  \begin{bmatrix} \text{---} & \mathbf{\Psi}^T(\y_1) & \text{---} \\  & ... &  \\ \text{---} & \mathbf{\Psi}^T(\y_M) & \text{---} \end{bmatrix}.
    \ee
    Let $\mathbf{B}=\mathbf{\Psi}_p$, and $\b=\mathbf{0}$. 
    Then the minimizer of \eqref{equ:primalmin} would be the coefficients of an approximation polynomial which is non-negative at points $\y_1$,...,$\y_M$. In this case, the number of constraints $C=M$.


\section{Efficient Iterative Methods}
\label{sec:method}
We now discuss splitting algorithms for solving convex optimization \eqref{equ:primalmin}. In particular, we focus on the Fast Iterative Shrinkage-Thresholding Algorithm (FISTA), \RV{which is an efficient simple algorithm if applied to a proper problem set up.}

\subsection{The primal and dual problems}
Recall that the primal problem \eqref{equ:primalmin} is given by
\be 
   \tag{\textbf{P}}
   \label{primal}
   \min_{\mathbf{c} \in \mathbb{R}^N} ~~\frac{1}{2}\alpha\norm{\mathbf{\Psi}_a \mathbf{c}-\mathbf{f}}_2^2 + \iota_{\Lambda}( \mathbf{B}\mathbf{c} - \b), 
\ee
Then \eqref{primal} can be written as 
\[ \min_{\mathbf{c} \in \mathbb{R}^N} g(\mathbf{c})+h (\mathbf{B}\mathbf{c} - \b),\]
where the two functions $g:\mathbb{R}^N \mapsto \mathbb{R}$ and $h:\mathbb{R}^C \mapsto \mathbb{R}$ are
\be
    \label{equ:not}
    g(\mathbf{c}) = \frac{1}{2}\alpha\norm{\mathbf{\Psi}_a \mathbf{c}-\mathbf{f}}_2^2,\qquad h(\x)=\iota_{\Lambda}( \x).
\ee
For any convex function $h$, its convex conjugate is defined as
\[h^\star(\y):=\sup_{\x \in \mathbb{R}^C} \left\{\left\langle \y, \x\right\rangle-h(\x)\right\}.\]    
The primal problem \eqref{primal} is equivalent to the following primal-dual form \eqref{equ:probpd} and dual form \eqref{equ:probd} (e.g., see \cite{rockafellar1967duality}):
\begin{align}
  \min_{\mathbf{c}\in \mathbb{R}^N}\max_{\mathbf{u}\in \mathbb{R}^C}& ~~\left[ \langle\mathbf{B} \mathbf{c}-\b,\mathbf{u}\rangle-h^\star(\mathbf{u})+g(\mathbf{c}) \right] \tag{\textbf{P-D}}\label{equ:probpd},\\
  -\min_{\mathbf{u}\in \mathbb{R}^C}& ~~g^\star(-\mathbf{B}^T\mathbf{u})+\b^T\mathbf{u}+h^\star(\mathbf{u}), \tag{\textbf{D}}\label{equ:probd}
\end{align}
where $\mathbf{u}\in \mathbb{R}^C$ is the dual variable.
For the two functions in \eqref{equ:not},   the conjugate functions $g^\star:\mathbb{R}^N \mapsto \mathbb{R}$, $h^\star:\mathbb{R}^C \mapsto \mathbb{R}$ are given by
\begin{align}
  g^\star(\mathbf{c}) &=\frac{1}{2\alpha}\left(\z+\mathbf{c}\right)^T\K^\dagger\left(\z+\mathbf{c}\right)-\frac{1}{2}\alpha\f^T\f\label{equ:cgg},\\
  h^\star(\mathbf{x}) &=\iota_{\Lambda^\star}(\mathbf{x}),~~~\text{with }\Lambda^\star = \{ \mathbf{\x}: \mathbf{x}\leq \mathbf{0}\}, \label{equ:cgh}
\end{align}
where $\z=\alpha\mathbf{\Psi}_a^T \f$, $\K=\mathbf{\Psi}_a^T\mathbf{\Psi}_a$, 
 and $\K^\dagger$ is the Moore-Penrose pseudo inverse $\K$. Let
\[ G^*(\mathbf{u}) = g^\star(-\mathbf{B}^T\mathbf{u})+\b^T\mathbf{u},\]   
then the dual problem \eqref{equ:probd} can be written as:
\be
   -\min_{\mathbf{u}\in \mathbb{R}^C} ~~G^*(\mathbf{u})+h^\star(\mathbf{u}), \label{equ:probd2} 
\ee
where 
\be
   G^*(\mathbf{u}) = \frac{1}{2\alpha}\left(\z-\mathbf{B}^T\mathbf{u}\right)^T\K^\dagger\left(\z-\mathbf{B}^T\mathbf{u}\right)-\frac{1}{2}\alpha\f^T\f+\b^T\mathbf{u}. \label{equ:Fu} 
\ee

\subsection{Closed convex proper functions}

\RV{In this subsection, we explain why the indicator function of a set must be defined using $+\infty$  in \eqref{indicator}.
For a set $\Omega\subset \mathbb R^n$, consider a function $$I_\Omega(\x)=\begin{cases}
0, & \x \in \Omega, \\ 
M, & \x \notin \Omega.
\end{cases}$$ for a very large number $M$. Then $I_\Omega(\x)$ is a well defined function on the whole space $\mathbb R^n$, and it may seem that $I_\Omega(\x)$ can serve the same purpose numerically as the indicator function \eqref{indicator}. However, a convex function well defined on the whole space $\mathbb R^n$ must be a continuous function \cite[Corollary 10.1.1]{rockafellar}, thus the function $I_\Omega(\x)$ cannot be convex on $\mathbb R^n$. To this end, one must consider a closed convex proper function, which will be defined as follows.

Notice that the indicator function defined in \eqref{indicator} should be regarded 
as an {\it extended} function $$f:\mathbb R^n \longrightarrow \mathbb R\cup \{\pm\infty\}.$$

For an extended function $f:\mathbb R^n \longrightarrow \mathbb R\cup \{\pm\infty\}$, its {\it epigraph} is defined as $$\mbox{epi} f=\{(\x,a)\in \mathbb R^n\times \mathbb R: f(\x)\leq a\}.$$
An extended function is called {\it closed} if its
epigraph is a closed set in $\mathbb R^{n+1}$. 
The {\it domain} of an extended function is denoted as $$\mbox{dom} f=\{\x\in \mathbb R^n: f(\x)\in  \mathbb R\}.$$
The indicator function of a set $\Omega$ defined in \eqref{indicator} is a closed extended function if and only if $\Omega$ is a closed set, see \cite{beck2017first}. 

A {\it convex} extended function is defined as an extended function satisfying:
\[ f(\lambda \x+(1-\lambda)\y )\leq \lambda f(\x)+(1-\lambda) f(\y),\quad \forall \x,\y \in \mbox{dom} f,~ \forall \lambda \in (0,1).\]

An extended function is called {\it proper} if it never maps negative infinity.
For instance, the  indicator function  \eqref{indicator} is a proper function. 
Thus if $\Omega$ is a closed convex set, then the indicator function \eqref{indicator} is a 
closed convex proper function, to which many results about convex functions well defined on the whole space $\mathbb R^n$ can be extended.

It can be proven that a closed extended function is also a lower semi-continuous function \cite[Theorem 2.6]{beck2017first}. Thus for a closed convex set $\Omega$, the indicator function \eqref{indicator} is also a 
lower semi-continuous convex proper function.

}

\subsection{First order splitting methods}
\label{sec:FISTA}
Both the primal problem \eqref{primal} and the dual problem  in the form of \eqref{equ:probd2} can be written as a general composite minimization 
\be \label{equ:gminprob}
    \min_{\x} ~~g(\x) + h(\x),
\ee
where the functions $g(\x)$ and $h(\x)$ are \RV{convex closed proper functions.} 
\RV{For a convex closed proper function  $g(\x)$ such as   the indicator function \eqref{indicator},
its subdifferential $\partial g(\x_0)$ at a point $\x_0$ is defined as a set of slopes of subtangent lines:
\[ \partial g(\x_0)=\{ \v\in\mathbb R^n: g(\x)\geq g(\x_0)+\langle \v, \x-\x_0\rangle,\quad \forall \x \in \mbox{dom} (g)\}. \]
An element in the subdifferential set $\partial g$ is called a subgradient.
If $g(\x)$ is differentiable at $\x_0$, then the subgradients coincide with the gradient, i.e., $\partial g(\x_0)=\{ \nabla g(\x_0)\}$. For example,  for $f(x)=|x|$, the subdifferential set  is 
$ \partial f(x)=\begin{cases}
    \{1 \}, & x>0\\
    \{-1 \}, & x<0\\
    [-1,1], & x=0
\end{cases}$.
}

Let $\partial g$ and $\partial h$ be their subdifferentials,
$I$ the identity operator. The proximal operators of these two functions,
i.e.,  the resolvents of subdifferentials, are
\begin{align}
   \label{equ:prox}
   \prox_{g}^{\gamma} (\x) =(I +\gamma \partial g)^{-1}(\x) = \operatorname{argmin}_\z \gamma g(\z)+\frac{1}{2}\|\z-\x\|_2^2, \quad \gamma>0,\\
   \prox_{h}^{\gamma} (\x) =(I +\gamma \partial h)^{-1}(\x) = \operatorname{argmin}_\z \gamma h(\z)+\frac{1}{2}\|\z-\x\|_2^2, \quad \gamma>0.
\end{align}

\RV{For a convex closed proper function $g(\x)$, $ \gamma g(\z)+\frac{1}{2}\|\z-\x\|_2^2$ is a closed strongly convex proper function and it has a unique minimizer \cite[Theorem 5.25]{beck2017first}, thus the  proximal operator   $ \prox_{g}^{\gamma} (\x) $ is a well defined operator.}

Assume the proximal operators have explicit formulae or can be efficiently approximated.
\RV{If $g(x)$ is differentiable but $h(x)$ is not differentiable, then the simplest method for \eqref{equ:gminprob} is the subgradient method:
\[ \x_{k+1}=\x_k-\gamma_k (\nabla g(\x_{k})+\v_k),\quad    \v_k\in \partial h(\x_k), \]
where $\gamma_k$ is a step size and any subgradient $\v_k$ can be chosen. However, the subgradient method may converge very slowly. Instead, a method using the proximal operator $\prox_{h}^{\gamma} (\x)$ is usually much faster.
The simplest splitting method to use the proximal operator is the forward-backward splitting:
}
\begin{equation}
 \x_{k+1}=(I +\gamma \partial h)^{-1}(I -\gamma \nabla g)(\x_{k})=\prox_{h}^{\gamma}(\x_{k}-\gamma \nabla g(\x_{k})).
 \label{F-B-splitting}
\end{equation} 
The splitting method \eqref{F-B-splitting} is also referred to as the {\it proximal gradient method}, or {\it projected gradient method} when $\prox_{h}^\gamma$ is a projection operator. The {\it fast proximal gradient method} \cite{nesterov2013gradient} is also called
 the Fast Iterative Shrinkage-Thresholding Algorithm (FISTA) \cite{beck2009fast}. For \eqref{equ:gminprob}, FISTA is given by
\be
\label{equ:FISTA}
\left\{\begin{array}{l}
\vspace{1mm} 
\x_{k+1}= \prox_{h}^\gamma (\y_{k} - \gamma \nabla g(\y_{k})), \\
\vspace{1mm}
t_{k+1} = \frac{1+\sqrt{1+4t_k^2}}{2}, \\
\vspace{1mm}
\y_{k+1} = \x_{k+1} + \left( \frac{t_k-1}{t_{k+1}} \right)(\x_{k+1}-\x_{k}),
\end{array}\right.
\ee
where $t_0=1$, $\x_0=\y_0$,  and $\gamma>0$ is a step size. The convergence rate of \eqref{F-B-splitting} can be proven $\mathcal O(k^{-1})$, and the convergence rate of FISTA \eqref{equ:FISTA} is $\mathcal O(k^{-2})$, when the step size is taken as $\gamma=\frac{1}{L}$, assuming that $\nabla g$ is Lipschitz continuous with Lipschitz constant $L$.

When both the proximal operators \st{are available} \RV{have explicit formulae or can be efficiently approximated}, we can also consider the splitting methods with two proximal operators. Define reflection operators as 
$$\mathtt{R}_{g}^\gamma=2\prox_{g}^\gamma-I,\quad \mathtt{R}_{h}^\gamma=2\prox_{h}^\gamma-I, $$
then the Douglas-Rachford splitting \cite{lions1979splitting, eckstein1992douglas} can be given as 
\begin{equation}
    \left\{\begin{array}{l}
        \vspace{1mm}
        \mathbf{y}_{k+1}= \lambda\frac{\mathtt{R}_g^\gamma\mathtt{R}_{h}^\gamma +\mathtt{I}}{2} \mathbf{y}_k+(1-\lambda)\mathbf{y}_k, \\
        \vspace{2mm}
        \mathbf{x}_{k+1}= \prox_{h}^\gamma (\mathbf{y}_{k+1}),
        \end{array}\right.
        \label{DR-splitting}
\end{equation}
where $\lambda\in(0,2]$ is a relaxation parameter. When $\lambda=2$, \eqref{DR-splitting} converges only if at least one of the two functions $g, h$ is strongly convex. The Douglas-Rachford splitting \eqref{DR-splitting}  is equivalent to the popular alternating direction method of multipliers (ADMM) method  \cite{gabay83} and the split Bregman method \cite{Goldstein:2009:SBM:1658384.1658386} on the equivalent Fenchel dual problem of \eqref{equ:gminprob} if using special step sizes,  see \cite{demanet2016eventual} and the references therein. When using ADMM and the split Bregman on the primal problem \eqref{equ:gminprob}, it is equivalent to using Douglas-Rachford splitting on the Fenchel dual problem of \eqref{equ:gminprob}. Usually, there is no significant difference in numerical performance between using the same splitting method on the primal problem and the dual problem.  

\begin{table}[htb!]
\centering
$\begin{array}{cl}
\toprule
 \text{\textbf{Method}} & 
\multicolumn{1}{c}{\hspace{-12mm}\text{\textbf{Iteration Schemes}}} \\
\midrule
\addlinespace
\text{Projected Gradient for \eqref{equ:probd2}} &  
        \left\{\begin{array}{l}
        \vspace{1mm}
        \mathbf{u}_{k}= \mathbf{u}_{k} - \eta \nabla G^*(\mathbf{u}_{k}), \\
        \vspace{2mm}
        \mathbf{u}_{k}= \prox_{h^\star}^\eta (\mathbf{u}_{k}).
        \end{array}\right.\\
\addlinespace
\midrule
\addlinespace
\text{Douglas-Rachford for \eqref{equ:probd2}}  &  
        \left\{\begin{array}{l}
        \vspace{1mm}
        \mathbf{p}_{k+1}= \frac{\mathtt{R}_{G^*}^\gamma\mathtt{R}_{h^\star}^\gamma +\mathtt{I}}{2} \mathbf{p}_k, \\
        \vspace{2mm}
        \mathbf{u}_{k+1}= \prox_{h^\star}^\gamma (\mathbf{p}_{k+1}).
        \end{array}\right.\\
\addlinespace
\midrule
\addlinespace
\text{Fast PDHG for \eqref{equ:probpd}}  & 
    \left\{\begin{array}{l}
    \vspace{1mm}
    \mathbf{u}_{k+1}=\prox_{h^\star}^{\tau_k}\left(\mathbf{u}_k+\tau_k (\mathbf{B} \bar{\mathbf{c}}_k-\b)\right), \\
    \vspace{1mm}
    \mathbf{c}_{k+1}=\prox_{g}^{\eta_k}\left(\mathbf{c}_k-\eta_k \mathbf{B}^T \mathbf{u}_{k+1}\right), \\
    \vspace{1mm}
    \theta_k=1 / \sqrt{1+2 \mu \eta_k}, \eta_{k+1}=\theta_k \eta_k, \tau_{k+1}=\tau_k / \theta_k, \\
    \bar{ \mathbf{c}}_{k+1}= \mathbf{c}_{k+1}+\theta_k\left( \mathbf{c}_{k+1}- \mathbf{c}_k\right).
    \end{array}\right.\\
\addlinespace
\bottomrule
\end{array}$ 
\caption{Examples of applying popular first-order splitting algorithms for solving \eqref{primal}, for which $\mathbf{c}$ is the primal variable, $\mathbf{u}$ is the dual variable, $\mathbf{p}$ is an auxiliary variable in Douglas-Rachford splitting and $\mu>0$ is a parameter for accelerated PDHG method. }
\label{tbl:scheme}
\end{table}

Another popular splitting method for using two proximal operators is the accelerated Primal-Dual Hybrid Gradient (PDHG) method \cite{chambolle2011first,pock2011diagonal}, for solving the primal-dual form \eqref{equ:probpd}.

 See \cite{beck2017first, ryu2022large} for a comprehensive introduction of these first-order algorithms.
In Table \ref{tbl:scheme}, we list several popular first-order algorithms that are used in our numerical tests to solve the optimization problem \eqref{equ:probd} (or in the form of \eqref{equ:probd2}) and \eqref{equ:probpd}.

\subsection{FISTA with restart on the dual}

 It is however inefficient to directly apply \eqref{equ:FISTA} to the problem \eqref{equ:primalmin}. The main technical difficulty is the computation cost of the proximal operator for $\iota_{\Lambda}( \mathbf{B}\mathbf{c} - \b)$ as a function of $\mathbf{c}$. The proximal operator of an indicator function is the projection operator to the domain it indicates. The domain $\{\mathbf{c}:\mathbf{B}\mathbf{c} - \b\geq 0\}$ is usually a high-dimensional hyper polygon defined by the linear inequalities. Although methods such as Dykstra's algorithm \cite{dijkstra2022note} can be used, the additional computational cost is undesired. We refer readers to \cite{zala2020structure} for a projection-based method for structure-preserving polynomial approximations. To circumvent the difficulty, we take advantage of the linear constraints and \RV{consider the equivalent dual problem \eqref{equ:probd2}}. This allows us to compute the proximal operator of a simple function with a closed-form formula for each step.

\begin{remark}
    Here we remark that the primal form \eqref{equ:primalmin} is in the space of $\mathbb{R}^N$, while the dual form is in $\mathbb{R}^C$. In practice, especially in high dimensions, we expect $N \gg C$. Thus the dual problem could significantly reduce the computational cost. For example, for a non-negative polynomial approximation, when $d$ is large, the cardinality of polynomial space $N$ can be exceedingly large even when polynomial order $n$ is very small. However, in practice, usually only a small set of discrete samples is needed to enforce the positivity, e.g., positivity is needed only at certain locations of interest, thus $C\ll N$.
\end{remark}

Since the function $G^*(\mathbf{u})$ in \eqref{equ:probd2} is quadratic with respect to $\mathbf{u}$\RV{, the} derivatives and proximal operators can be written explicitly as:
\begin{align}
    \nabla G^*(\mathbf{u}) & = \frac{1}{\alpha} \left( \mathbf{B}\K^\dagger\mathbf{B}^T \mathbf{u} -\mathbf{B}\K^\dagger\z \right) + \b, \label{equ:dF} \\
    \prox_{G^*}^\gamma (\mathbf{u}) & = \left(I+\frac{\gamma}{\alpha}\mathbf{B}\K^\dagger\mathbf{B}^T\right)^{-1}\left[ \mathbf{u}+\gamma\left( \frac{1}{\alpha}\mathbf{B}\K^\dagger\z-\b \right) \right].\label{equ:pF}
\end{align}

Since $h^\star(\mathbf{u})$ in \eqref{equ:probd2} is an indicator function on the negative half-space of $\mathbb{R}^C$. Its proximal operator is simply the cut-off operator:
\begin{equation}
	\prox_{h^\star}^\tau (\mathbf{u}) = \left(\min\{u_i,0\}\right)_{i=1}^{C}.
\end{equation}

\RV{Now the FISTA method can be efficiently implemented on the dual problem \eqref{equ:probd2}.
  In practice, FISTA method can be further accelerated by various restarting strategies \cite{o2015adaptive, beck2017first}, which is accomplished by reiterating the sequence $t_k$ from the starting point $t_0$ after some iterations. The simplest restarting scheme is to repeat the standard FISTA algorithm \eqref{equ:FISTA}
with a fixed frequency. }

Here, we consider the following two criteria in \cite{o2015adaptive}:
\begin{align}
  & g(\x_k) > g(\x_{k-1}) \tag{A}\label{equ:restartcond1},\\
  & [\nabla g\left(\y_{k-1}\right)]^T\left(\x_k-\x_{k-1}\right)>0. \tag{B}\label{equ:restartcond2}
\end{align} 
The FISTA with adaptive restart, hereafter referred to as r-FISTA, usually achieves much faster convergence than the standard FISTA \cite{beck2017first}. In our numerical tests, the r-FISTA applied to the dual problem is superior to all the methods listed in Table \ref{tbl:scheme} using the same step size.

Then our main algorithm is to apply FISTA with adaptive restart (r-FISTA) on the dual problem \eqref{equ:probd} which results in an optimization in $\mathbb{R}^C$. After obtaining the minimizer $\mathbf{u}^\star$, we recover the primal minimizer $\mathbf{c}^\star$ by the primal-dual relation:
\be
    \K \mathbf{c}^\star = \frac{1}{\alpha }(\z- \mathbf{B}^T \mathbf{u}^\star).
\ee
We summarize the above algorithm in Algorithm \ref{alg:FISTAwRe}.

\begin{algorithm}
\caption{Solve \eqref{equ:probd} by FISTA with Adaptive Restart}
\label{alg:FISTAwRe}
\begin{algorithmic}[1]
\REQUIRE Matrix $\mathbf{\Psi}_a$, $\mathbf{B}$; Vectors $\mathbf{f}$, $\mathbf{b}$; Constants $\alpha$, $\eta$ (step-size). Functions $G^*(\mathbf{u})$, $\nabla G^*(\mathbf{u})$; Stopping criteria $S$. Use condition \eqref{equ:restartcond1} for restart.


\hspace*{-1.6\algorithmicindent}\textbf{Output:} The basis coefficients $\mathbf{c}$.
\vspace{1mm}

\STATE{Compute $\K=\mathbf{\Psi}_a^T\mathbf{\Psi}_a$, $\K^\dagger$, $\z=\alpha\mathbf{\Psi}_a^T \f$;}
\STATE{Initialize $\mathbf{u}_0=\mathbf{u}_1=\mathbf{p}_1$, $t_1=1$, $k=1$;}

\WHILE{($S=$ \textbf{FALSE})}
    \STATE{$\mathbf{u}_{k+1} = \min (\mathbf{p}_k-\eta \nabla G^*(\mathbf{p}_k), \mathbf{0})$;}
    \STATE{$t_{k+1} = \frac{1+\sqrt{1+4t_k^2}}{2}$;}
    \STATE{$\mathbf{p}_{k+1} = \mathbf{u}_k +  \frac{t_k-1}{t_{k+1}}  (\mathbf{u}_{k+1}-\mathbf{u}_{k})$;}
    \IF{$G^*(\mathbf{u}_{k+1}) > G^*(\mathbf{u}_k)$}
        \STATE{$t_{k+1} = 1$;}
    \ENDIF
    \STATE{$k=k+1$;}
\ENDWHILE

\STATE{$ \mathbf{c} = \K^\dagger(\z- \mathbf{B}^T \mathbf{u}_k)/\alpha$.}

\end{algorithmic}
\end{algorithm}

We briefly discuss the computational cost of the proposed method. 
As a preparation step, the matrix $\mathbf{B}\K^\dagger\mathbf{B}^T$ of size $C \times C$ and two vectors $\mathbf{B}(\K^\dagger)^T\z$, $\mathbf{b}$ of size $C\times 1$ need to be computed and stored. 
Within each iteration of Algorithm \ref{alg:FISTAwRe}, it costs $\mathcal{O}(C^2)$ flops. 
\RV{As to the convergence rate, we mention some classical results of  of FISTA for convex optimization in the next subsection}. 
In particular, for the problem \eqref{equ:probd2},  the convergence rate of first-order methods is related to the ratio between the Lipschitz constant of $G^*$ and the strong convexity parameter of the cost function. Since $G^*$ is quadratic, such a ratio is naturally related to 
the condition number of the matrix $\mathbf{B}\K^\dagger\mathbf{B}^T$.
It is also interesting to explore its performance dependency on the parameters, such as the number of constraints $C$, number of approximation points $K$, polynomial order $n$, and most importantly, the dimension $d$. However, it is difficult to conduct such analysis, as the properties the Vandermonde-like matrix $\mathbf{\Psi}_a$ depend critically on the geometric distribution of the samples $\x_1$,...,$\x_K$. Consequently, we rely on numerical testing in Section \ref{sec:numex}. Through numerical tests in both low and high dimensions, we discover that the number of iteration steps required for numerical convergence satisfies $k \sim \mathcal{O}(C)$. Under this condition, the proposed algorithm \ref{alg:FISTAwRe} has computational complexity $\mathcal{O}(C^3)$ and requires $\mathcal{O}(C^2)$ storage of real numbers. 

The cost of computing matrices $\mathbf{\Psi}_a$, $\mathbf{K}$, $\K^\dagger$ grows dramatically as dimension $d$ increases. This is due to the cardinality of the polynomial space $N \approx d^n/n!$ for large $d$. However, these computations need to be done only once and can be efficiently implemented in parallel, e.g.,  on modern Graphics Processing Units (GPUs).
    In our numerical tests in Section \ref{sec:numex}, the number of iterations needed for convergence does not increase when dimension $d$ increases.   When $N$ dramatically grows, the matrix $\mathbf{B}$ of the size of $C \times N$ becomes larger,
     but the size of $\mathbf{B}\K^\dagger\mathbf{B}^T$ does not grow much when we use a mild number of constraints $C$.

\subsection{The convergence rate of FISTA methods}

\RV{The standard FISTA method has the following provable $\mathcal{O}(1/k^2)$ convergence rate when applied to problem \eqref{equ:probd2}. 
By   \cite[Theorem 10.34]{beck2017first}, we have
\begin{theorem} \label{thm:FISTAn}
Let $\left\{\mathbf{u}_k\right\}_{k \geq 0}$ be the sequence generated by FISTA \eqref{equ:FISTA} for solving problem \eqref{equ:probd2}. Then for any $k \geq 1$,
$$
F\left(\mathbf{u}^k\right)-F_{\mathrm{opt}} \leq \frac{2 L\left\|\mathbf{u}^0-\mathbf{u}^*\right\|^2}{(k+1)^2},
$$
where $L$ is the maximum eigenvalue of the matrix $\mathbf{B}\K^\dagger\mathbf{B}^T$, $F(\mathbf{u}) := G^*(\mathbf{u}) + h^\star(\mathbf{u})$ is the objective function defined in \eqref{equ:probd2} with optimal solution $\mathbf{u}^*$ and optimum $F_{\mathrm{opt}}$.
\end{theorem}

 FISTA with a fixed restarting frequency is described  in Algorithm \ref{alg:FISTAwRe-2}.

\begin{algorithm}
\caption{Solve \eqref{equ:probd} by FISTA with  a fixed restarting frequency $N$}
\label{alg:FISTAwRe-2}
\begin{algorithmic}[1]
\REQUIRE Matrix $\mathbf{\Psi}_a$, $\mathbf{B}$; Vectors $\mathbf{f}$, $\mathbf{b}$; Constants $\alpha$, $\eta$ (step-size). Functions $G^*(\mathbf{u})$, $\nabla G^*(\mathbf{u})$; a fixed restarting frequency $N$.


\hspace*{-1.6\algorithmicindent}\textbf{Output:} The basis coefficients $\mathbf{c}$.
\vspace{1mm}

\STATE{Compute $\K=\mathbf{\Psi}_a^T\mathbf{\Psi}_a$, $\K^\dagger$, $\z=\alpha\mathbf{\Psi}_a^T \f$;}
\STATE{Initialize $\mathbf{u}_0=\mathbf{u}_1=\mathbf{p}_1$, $t_1=1$, $k=1$;}

\WHILE{($S=$ \textbf{FALSE})}
    \STATE{$\mathbf{u}_{k+1} = \min (\mathbf{p}_k-\eta \nabla G^*(\mathbf{p}_k), \mathbf{0})$;}
    \STATE{$t_{k+1} = \frac{1+\sqrt{1+4t_k^2}}{2}$;}
    \STATE{$\mathbf{p}_{k+1} = \mathbf{u}_k +  \frac{t_k-1}{t_{k+1}}  (\mathbf{u}_{k+1}-\mathbf{u}_{k})$;}
    \IF{$k \mod N=0$}
        \STATE{$t_{k+1} = 1$; }
        \STATE {$\mathbf{p}_{k+1}=\mathbf{u}_{k+1}$;}
    \ENDIF
    \STATE{$k=k+1$;}
\ENDWHILE

\STATE{$ \mathbf{c} = \K^\dagger(\z- \mathbf{B}^T \mathbf{u}_k)/\alpha$.}

\end{algorithmic}
\end{algorithm}

  With a stronger assumption on the  matrix $\mathbf{B}\K^\dagger\mathbf{B}^T$, FISTA with a fixed restarting frequency can be proven to converge linearly when applied to \eqref{equ:probd2}.
 By \cite[Theorem 10.41]{beck2017first}, we have
\begin{theorem}
Suppose the matrix $\mathbf{B}\K^\dagger\mathbf{B}^T$ is positive definite, and denote the minimum and maximum eigenvalue of $\mathbf{B}\K^\dagger\mathbf{B}^T$ by $\mu$ and $L$, condition number by $\kappa:=L/\mu$. Let $\{\mathbf{u}_k\}_{k\geq 0}$ be the sequence generated by the restarted FISTA method employed with a fixed restarting frequency $N=\lceil\sqrt{8 \kappa}-1\rceil$, then for any $k \geq 0$, the following convergence result holds
\begin{equation*}
    F\left(\mathbf{u}^k\right)-F_{\mathrm{opt}} \leq \frac{L \left\|\mathbf{u}^0-\mathbf{u}^*\right\|^2}{2}\left(\frac{1}{2}\right)^k,
\end{equation*}
where $F(\mathbf{u}) := G^*(\mathbf{u}) + h^\star(\mathbf{u})$ is the objective function defined in \eqref{equ:probd2} with optimal solution $\mathbf{u}^*$ and optimum $F_{\mathrm{opt}}$.
\end{theorem}

To implement this fixed-frequency restarting FISTA method, one has to choose the period that requires prior knowledge of condition number $\kappa$, which may not be available due to either lack of information or ill-conditioning of the system. The adaptive restart \cite{nesterov2013gradient, o2015adaptive} mitigates this issue by reiterating the sequence $t_k$ when certain criterion is met. 
 
\begin{remark} \label{rmk:vfista}
   When condition number $\kappa$ is known, one can improve standard FISTA by choosing an optimal stepsize instead of restarting. It is achieved by replacing $\frac{t_k-1}{t_{k+1}}$ by constant $\frac{\sqrt{\kappa}-1}{\sqrt{\kappa}+1}$ in \eqref{equ:FISTA}. This method, referred to as v-FISTA, can also achieve a provable linear convergence rate $\mathcal{O}((1-\frac{1}{\sqrt{\kappa}})^k)$ given a strongly-convex objective function \cite{beck2017first}. However, v-FISTA behaves worse than r-FISTA in our experiment due to the ill-conditioning system, which will be shown in numerical examples.
\end{remark}
}

\section{Numerical Examples}
\label{sec:numex}

In this section, we present several numerical examples to demonstrate the effectiveness of our proposed method. In all of our numerical tests, we use Legendre polynomials for the basis. The constants in Algorithm \ref{alg:FISTAwRe} are set to be: $\alpha=100$, $\eta=\alpha/|\sigma_{\mathtt{max}}(\mathbf{B}\K^\dagger\mathbf{B}^T)|$, where $\sigma_{\mathtt{max}}$ is the maximum eigenvalue. In our test, we use the following stopping criteria (i.e. $S$ in Algorithm \ref{alg:FISTAwRe}):
\begin{enumerate}
    \item Convergence of primal variable: $\norm{\mathbf{c}_k-\mathbf{c}_{k-1}}_2 \leq 1\times 10^{-14}$. Note that this is at the machine accuracy level.
    \item Satisfaction of constraints: $\mathbf{B}\mathbf{c}_k \geq \mathbf{b}$.
\end{enumerate} The dimensions of our examples include low dimensions $d=1$, $d=2$, intermediate dimension $d=10$, and high dimensions $d=100$, $d=200$.

We consider the following functions as testing. In one dimension ($d=1$), we consider the following $3$ functions:
\begin{align}
    & \textbf{Runge function}:  f_1(x)=\frac{101}{100}\left( \frac{1}{1+100x^2}-\frac{1}{101} \right); \label{equ:f1} \\
    & \textbf{Truncated sine function}:  f_2(x)=\left [\sin(\frac{\pi(x+1)}{2})-\sin(0.6\pi) \right ] \mathds{1}_{\{|x|<0.2\}}; \label{equ:tsin}\\
    & \textbf{Simple indicator function}:  f_3(x)=\mathds{1}_{\{x>0\}}, \label{equ:f3}
\end{align}
where $\mathds{1}_{\Omega}$ is the indicator function on set $\Omega$,
\[ \mathds{1}_{\Omega}(x)= \begin{cases}1, & x \in \Omega \\ 0, & x \notin \Omega\end{cases}.\]
In multi-dimension ($d>1$), we test the following 3 functions (\cite{genz1984testing}) that have been widely used for multi-dimensional function approximation:
\begin{align}
    & \hspace{-4mm}\textbf{Gaussian peak function}: f_4(\mathbf{x})=\exp \left(-\sum_{i=1}^d \sigma_i^2\left(\frac{x_i+1}{2}-\omega_i\right)^2\right); \label{equ:gaussian} \\
    & \hspace{-4mm}\textbf{Continuous peak function}: f_5(\mathbf{x})=\exp \left(-\sum_{i=1}^d \sigma_i\left|\frac{x_i+1}{2}-\omega_i\right|\right); \label{equ:f5}\\
    & \hspace{-4mm}\textbf{Corner peak function}: f_6(\mathbf{x})=\left(1+\sum_{i=1}^d \sigma_i \frac{\left(x_i+1\right)}{2}\right)^{-(d+1)}. \label{equ:f6}
\end{align} 

The parameters $\sigma_i$ and $\omega_i$, $i=1,...,d$ are designed to control the behavior of functions and specified in each example.

\subsection{One-dimensional Examples}
For one-dimensional cases, the domain of approximation is set to be $D=[-1,1]$. The samples for approximation $\x_i$ are chosen to be Chebyshev nodes with the number of samples $K=50$. The behavior of our proposed method is very similar in all the examples, so we only present a subset of the test results.

\subsubsection{Runge function}
\begin{figure}[htbp]
  \centering
  \label{fig:Runge}
  \includegraphics[width=.49\textwidth]{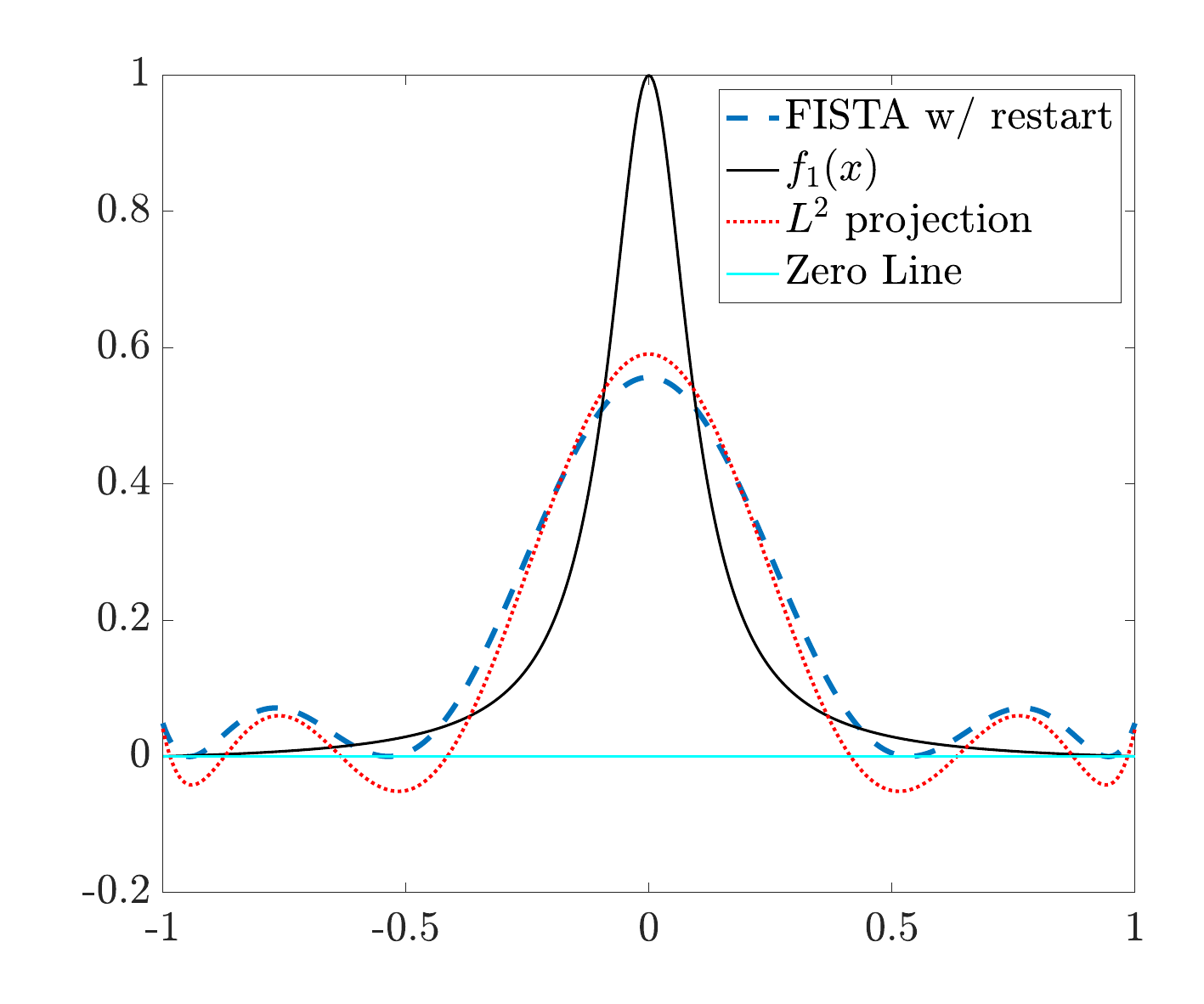}
  \includegraphics[width=.49\textwidth]{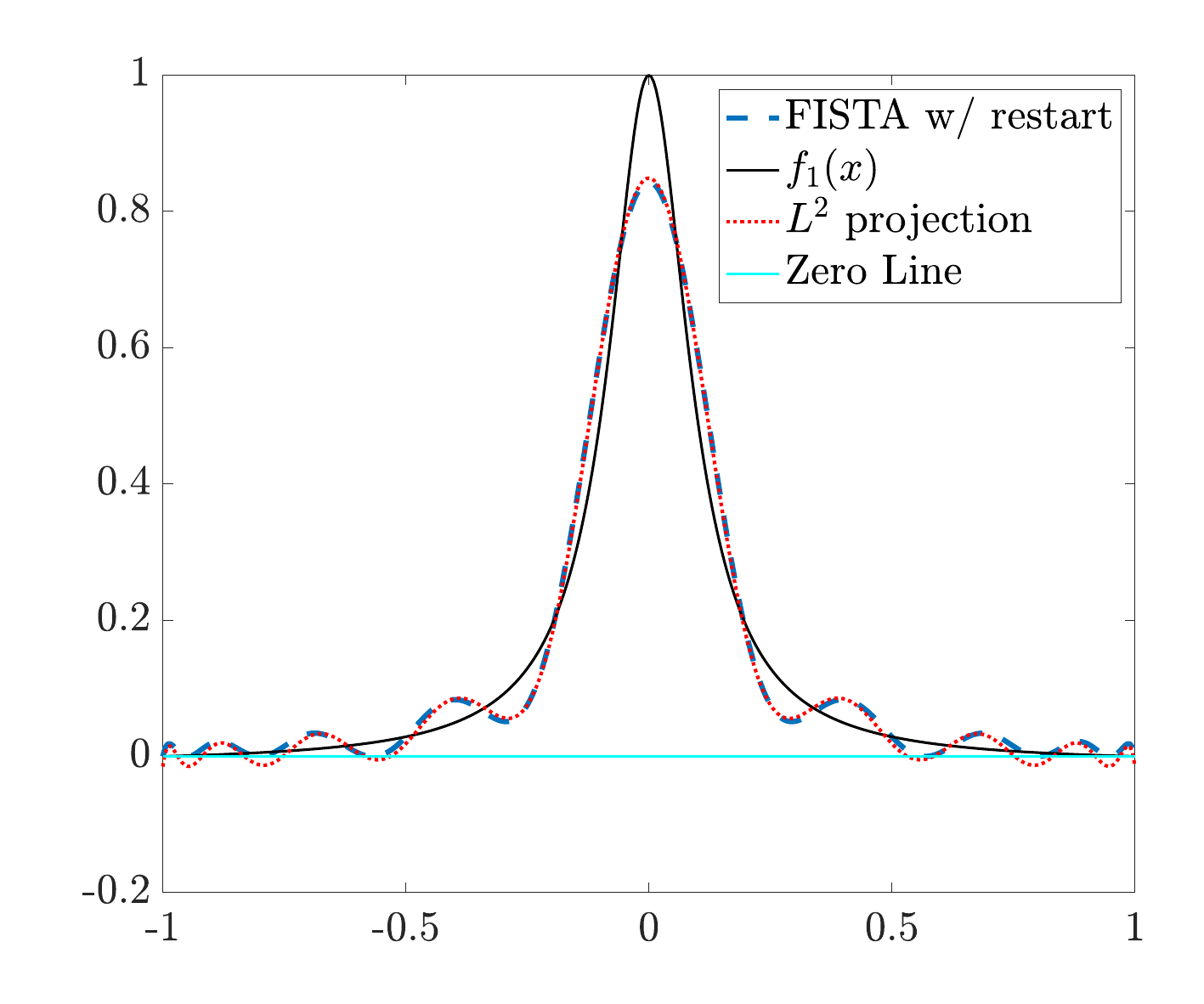}
  \caption{The approximated polynomials for $f_1(x)$ \eqref{equ:f1} with positivity constraints for $n=10$ (Left), $n=20$ (Right). The positive approximation is found via solving \eqref{primal} by restarted FISTA on \eqref{equ:probd2}.}
\end{figure}

\begin{figure}[htbp]
  \centering
  \label{fig:Runge_error}
  \includegraphics[width=.49\textwidth]{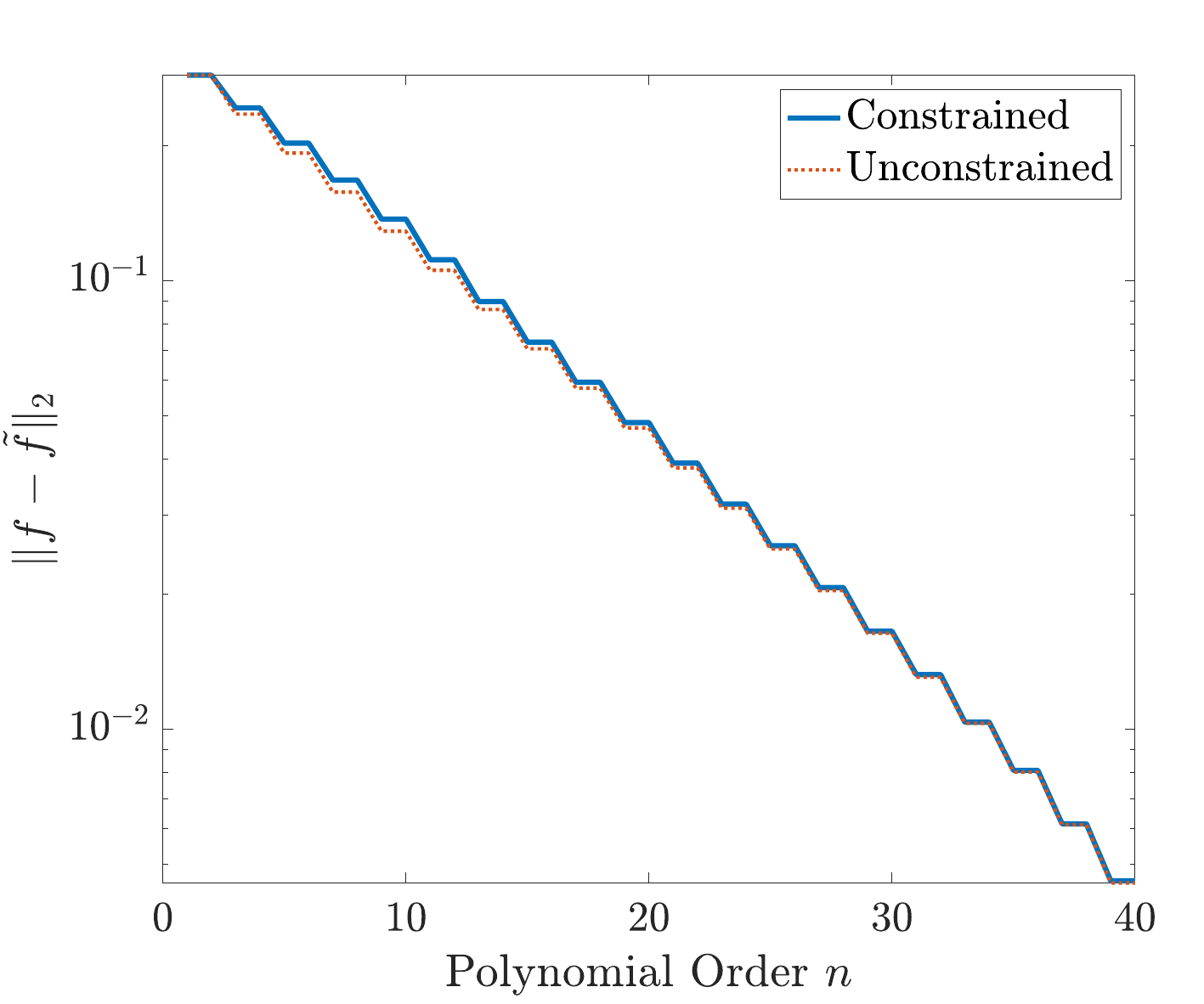}
  \includegraphics[width=.49\textwidth]{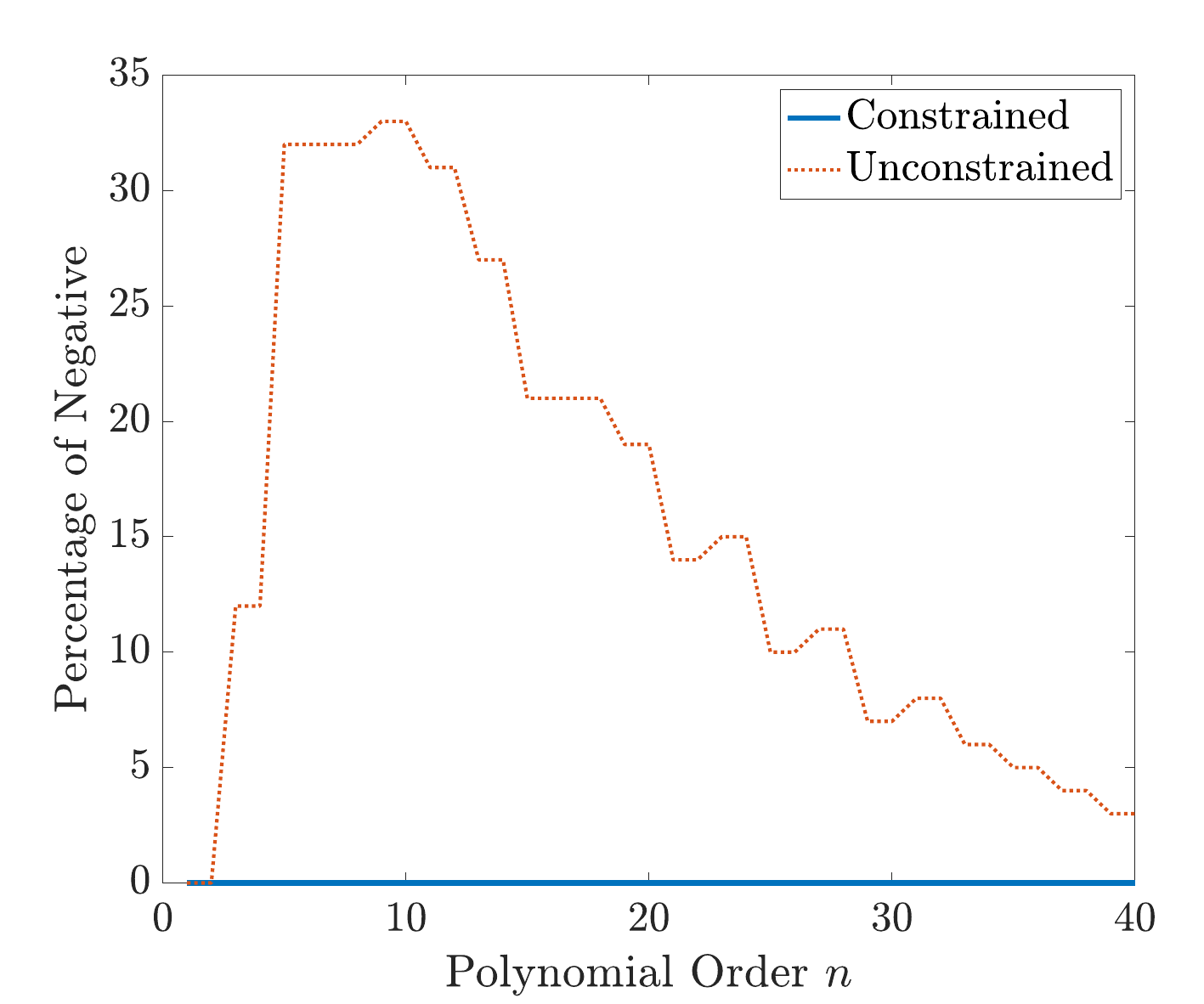}
  \caption{\RV{Comparasion for unconstrained approximation ($L^2$ projection) and non-negativity constrained approximation for $f_1(x)$ \eqref{equ:f1}. The convergence with respect to polynomial order (left) and percentage of negative points in $201$ non-negativity enforced sample points (right).}}
\end{figure}

\begin{figure}[htbp]
  \centering
  \label{fig:Runge_pts}
  \includegraphics[width=.32\textwidth]{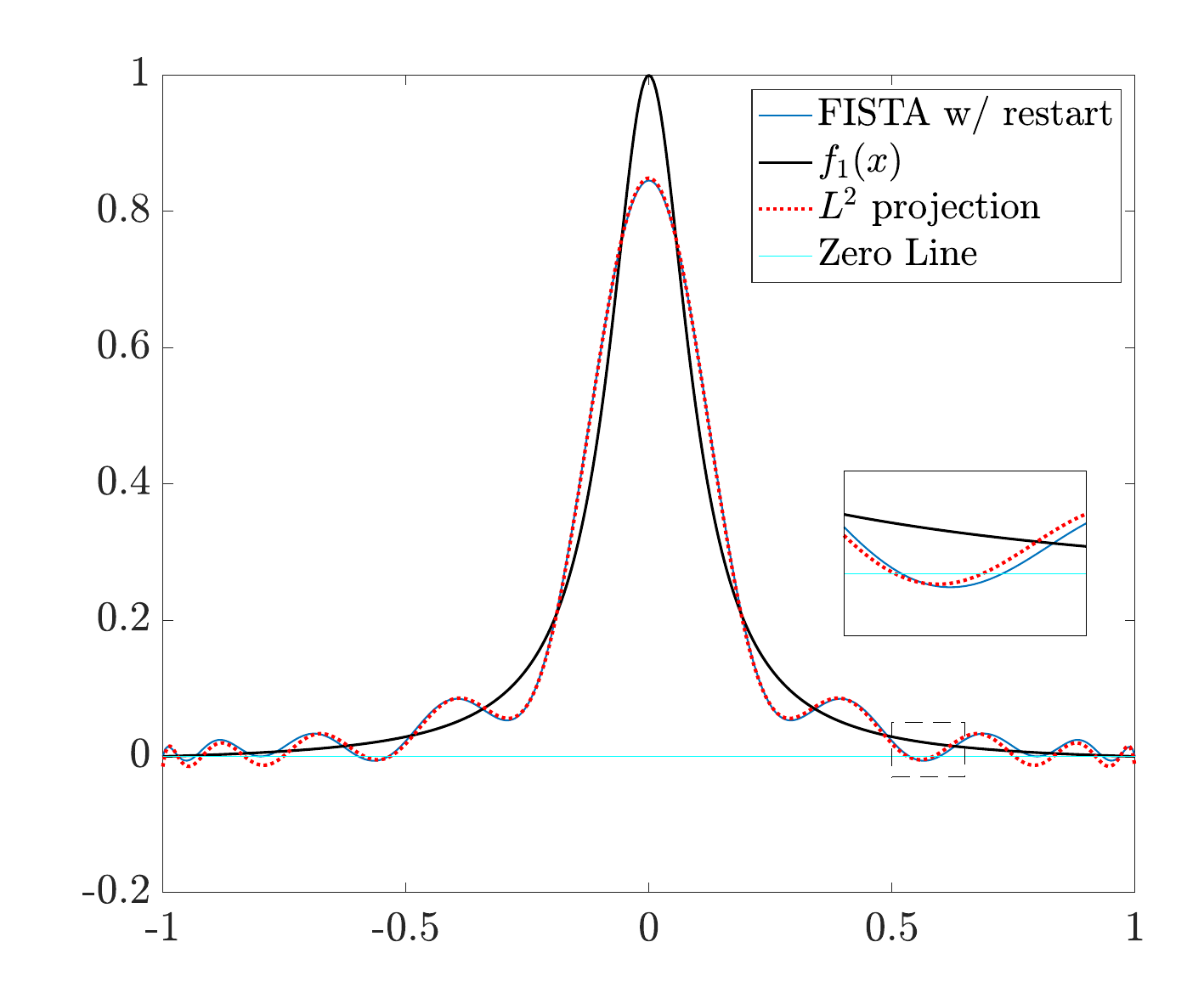}
  \includegraphics[width=.32\textwidth]{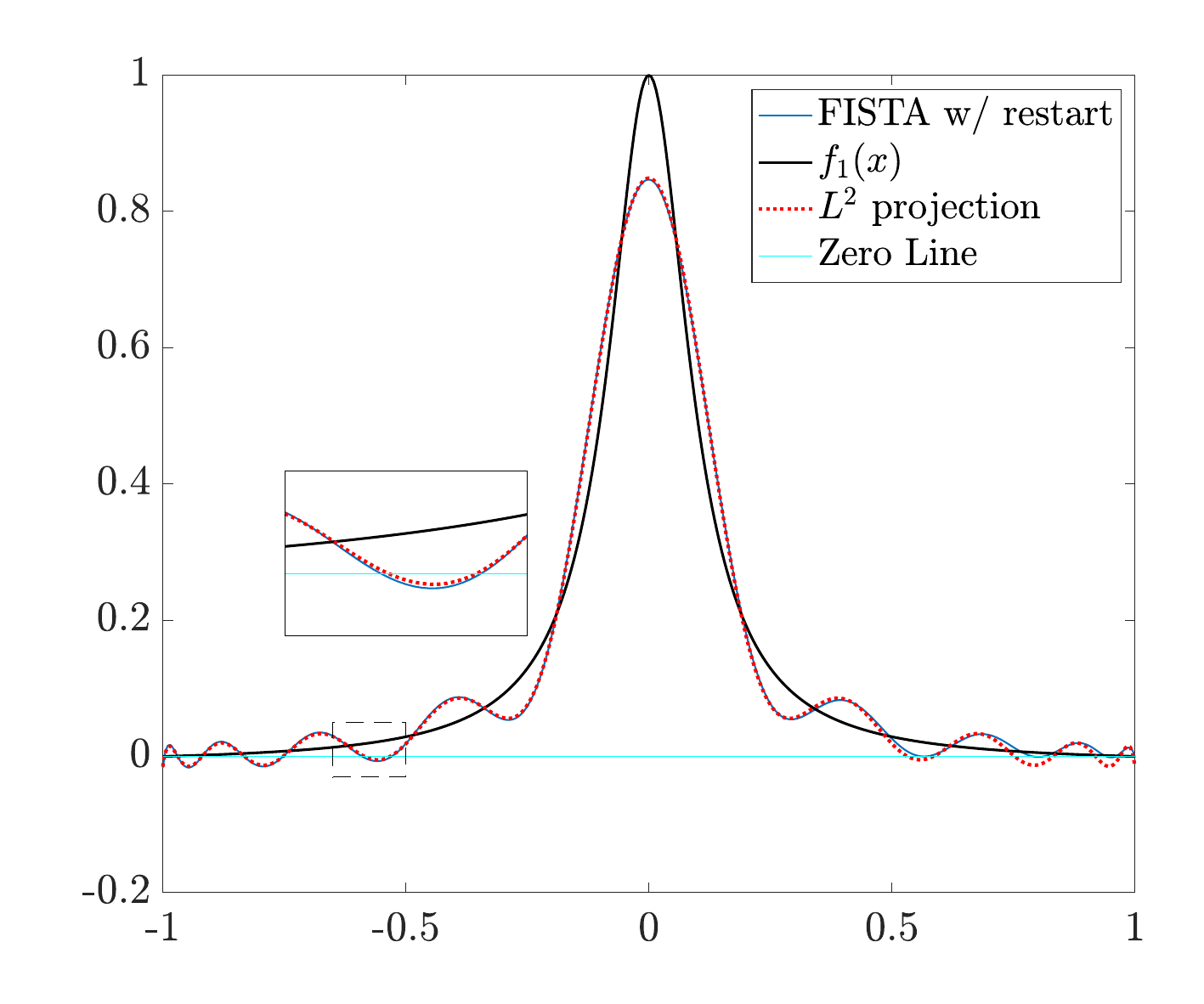}
  \includegraphics[width=.32\textwidth]{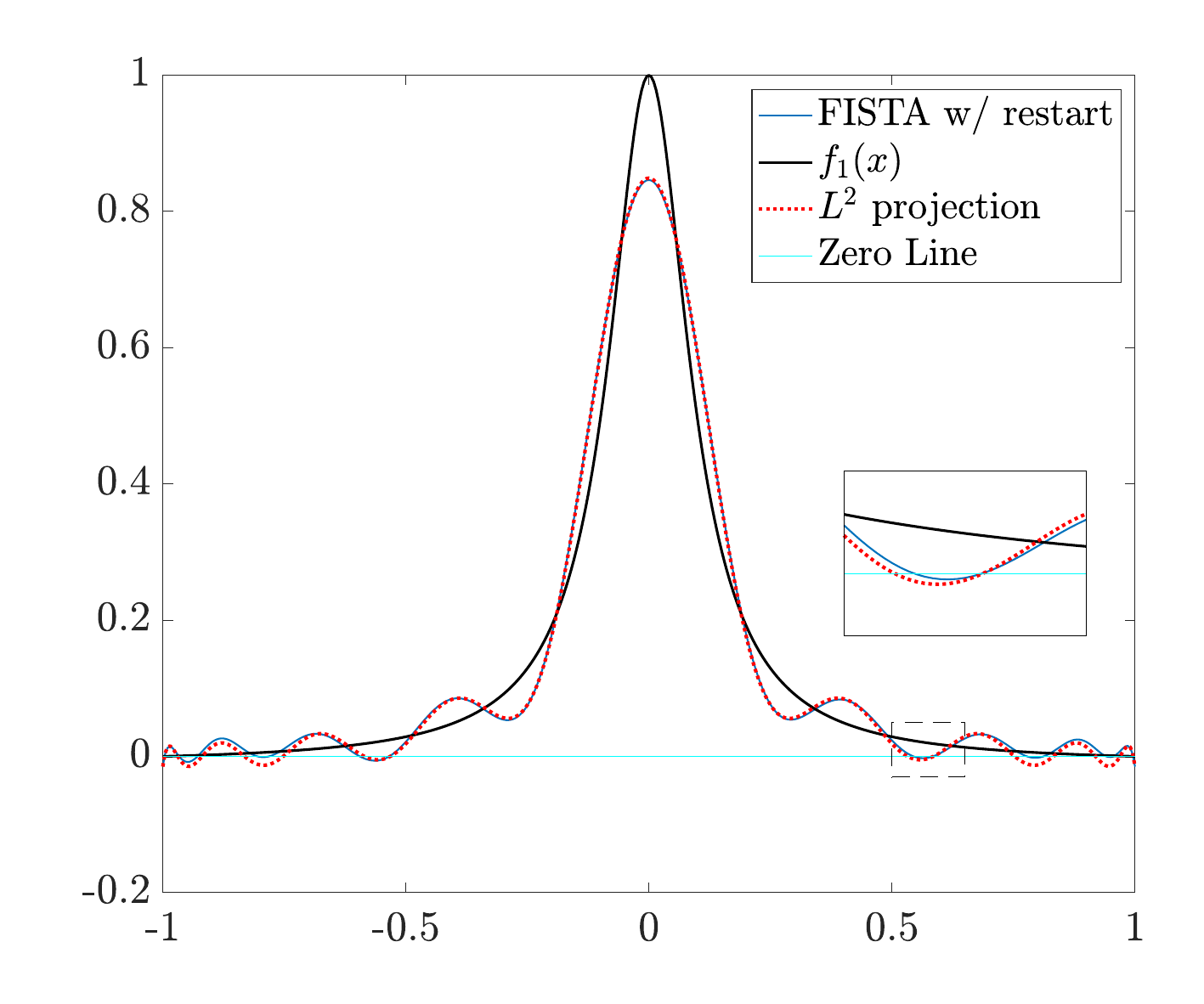}
  \caption{\RV{The approximated polynomials for $f_1(x)$ \eqref{equ:f1} with positivity constraints for $n=20$ with $M=31$ equispaced (left), Chebshev (middle) and random (right) constraint points. The areas enclosed by the dotted line are zoomed in for detailed presentation.}}
\end{figure}

We first consider the scaled Runge function $f_1(x)$ \cite{runge1901empirische}, which is positive. When the polynomial degree is large, its standard polynomial approximation produces oscillations around the edge, resulting in undesirable negative values. We take equidistant grid points, on which we enforce the positivity constraints, i.e. $\y_i=-1+\frac{2(i-1)}{M-1}$, $i=1,...,M$. Following the notation in Section \ref{sec:setup-constraint}, the matrix $\mathbf{B}=\mathbf{\Psi}_p$ is generated by the samples $\y_i$, $i=1,...,M$ and number of constraints $C=M$. For finding a positive approximation, we take $\mathbf{b}=\mathbf{0}+\epsilon$, with $\epsilon=10^{-5}$. The results for polynomial orders $n=10,20$ with $M=201$ are in Figure \ref{fig:Runge}. Compared with $L^2$ projection (red dot line), our approximation (blue dashed line) could almost preserve positivity on the domain. \RV{We test the convergence rate for approximation error $\|f-\tilde{f}\|_2$ with respect to polynomial order in Figure \ref{fig:Runge_error}. It is observed that both unconstrained estimation and constrained estimation converge at an exponential rate. The constrained error is slightly larger than the unconstrained error when the polynomial order is small, while this difference is eliminated as the polynomial order gets larger. The positivity constraints limit the optimality of approximation. As the polynomial degree increases, $L^2$ estimation suffers less from the positivity issue and becomes close to the constrained estimation. This conclusion is also reflected in the percentage of negativity points out of the $201$ constraint points, presented on the right of Figure \ref{fig:Runge_error}.

We also conduct tests for $f_1$ with different settings of nonnegativity-enforcing points. We test 3 sets of points: equidistant, Chebyshev, and uniformly distributed random points with $M=30$. The approximated polynomials are shown in Figure \ref{fig:Runge_pts}. It can be observed that the difference in the approximations is marginal. Though these approximated polynomials are guaranteed to be non-negative on the enforced points, it is still impossible to obtain non-negativity on all $x\in D$. This is also illustrated in Figure \ref{fig:Runge_pts}. In the figure, some violations can be seen in the area marked by dotted lines, which are zoomed in for clear presentation.
}

\subsubsection{Truncated sine function}
\begin{figure}[htbp]
  \centering
  \label{fig:Tsin}
  \includegraphics[width=.49\textwidth]{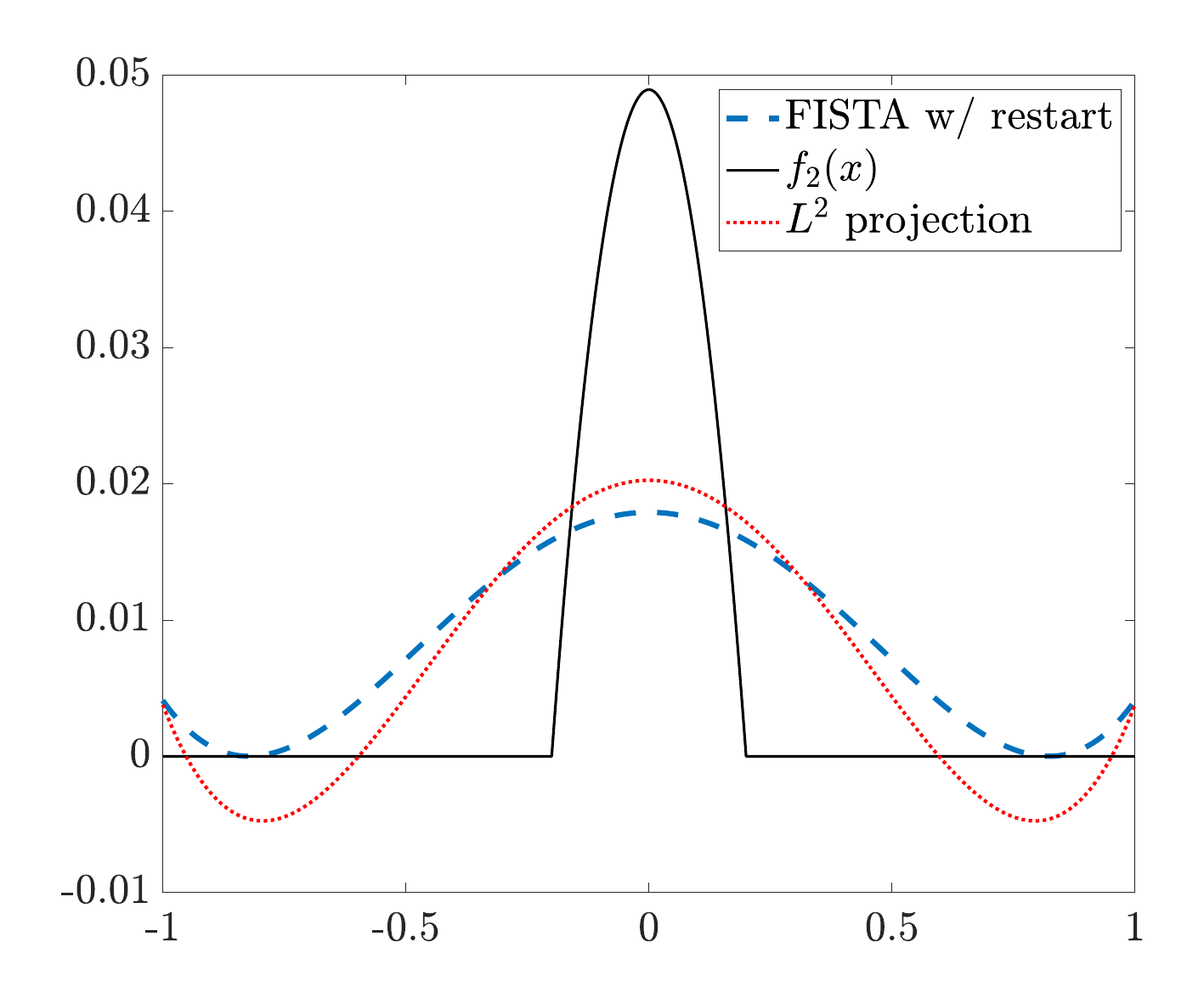}
  \includegraphics[width=.49\textwidth]{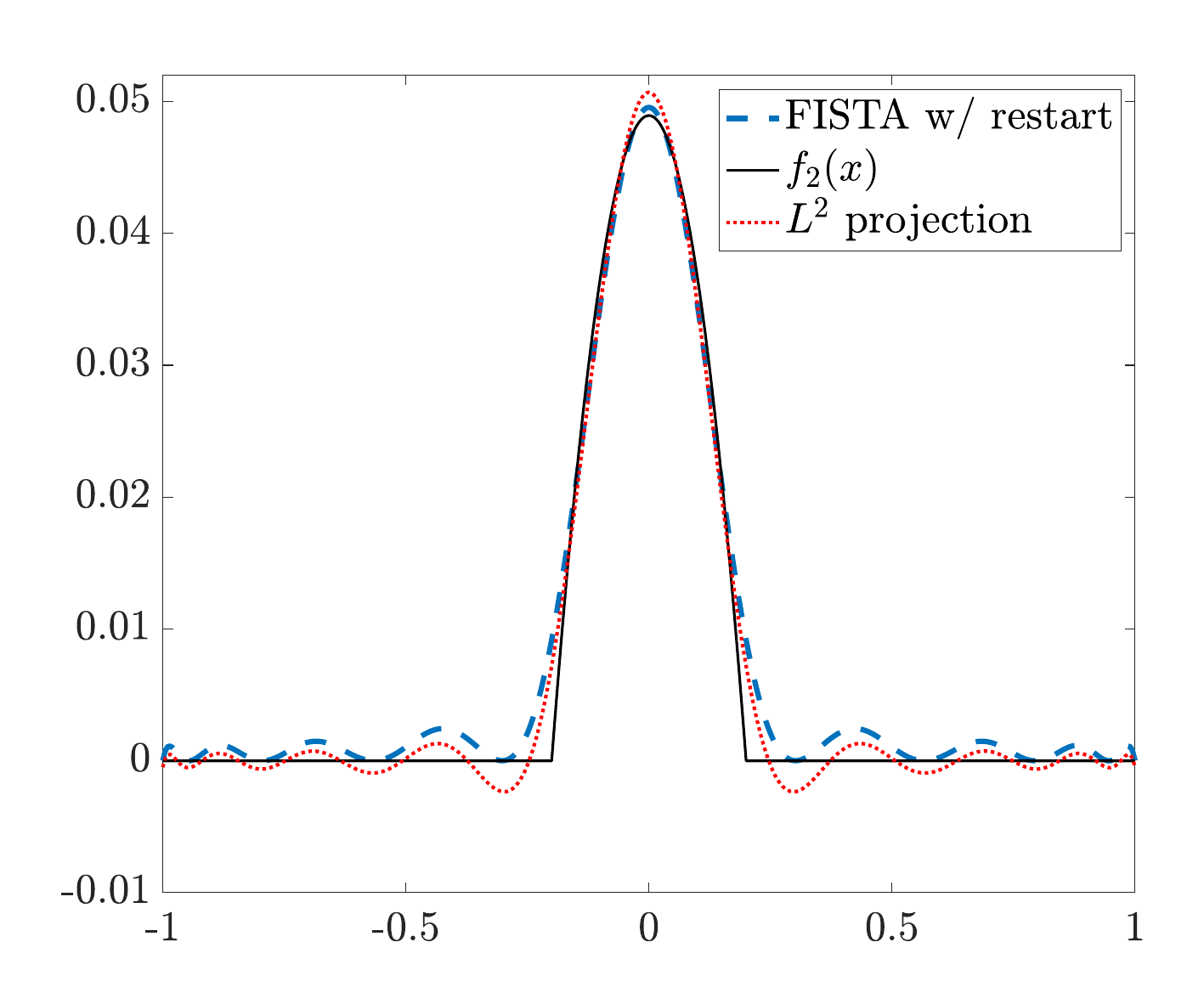}
  \caption{The approximated polynomials for $f_2(x)$ \eqref{equ:tsin} with positivity constraints for $n=5$ (Left), $n=20$ (Right). The positive approximation is found via solving \eqref{primal} by restarted FISTA on \eqref{equ:probd2}.}
\end{figure}

The truncated sine function $f_2(x)$ \eqref{equ:tsin} is a non-smooth, non-negative function in a bounded domain. We follow the same technical setting as $f_1(x)$ and seek a positive approximation. The results are shown in Figure \ref{fig:Tsin} for polynomial orders $n=5$ with $M=101$ and $n=20$ with $M=201$.

\begin{figure}[htbp]
  \centering
  \label{fig:Tsincompare}
  \includegraphics[width=.99\textwidth]{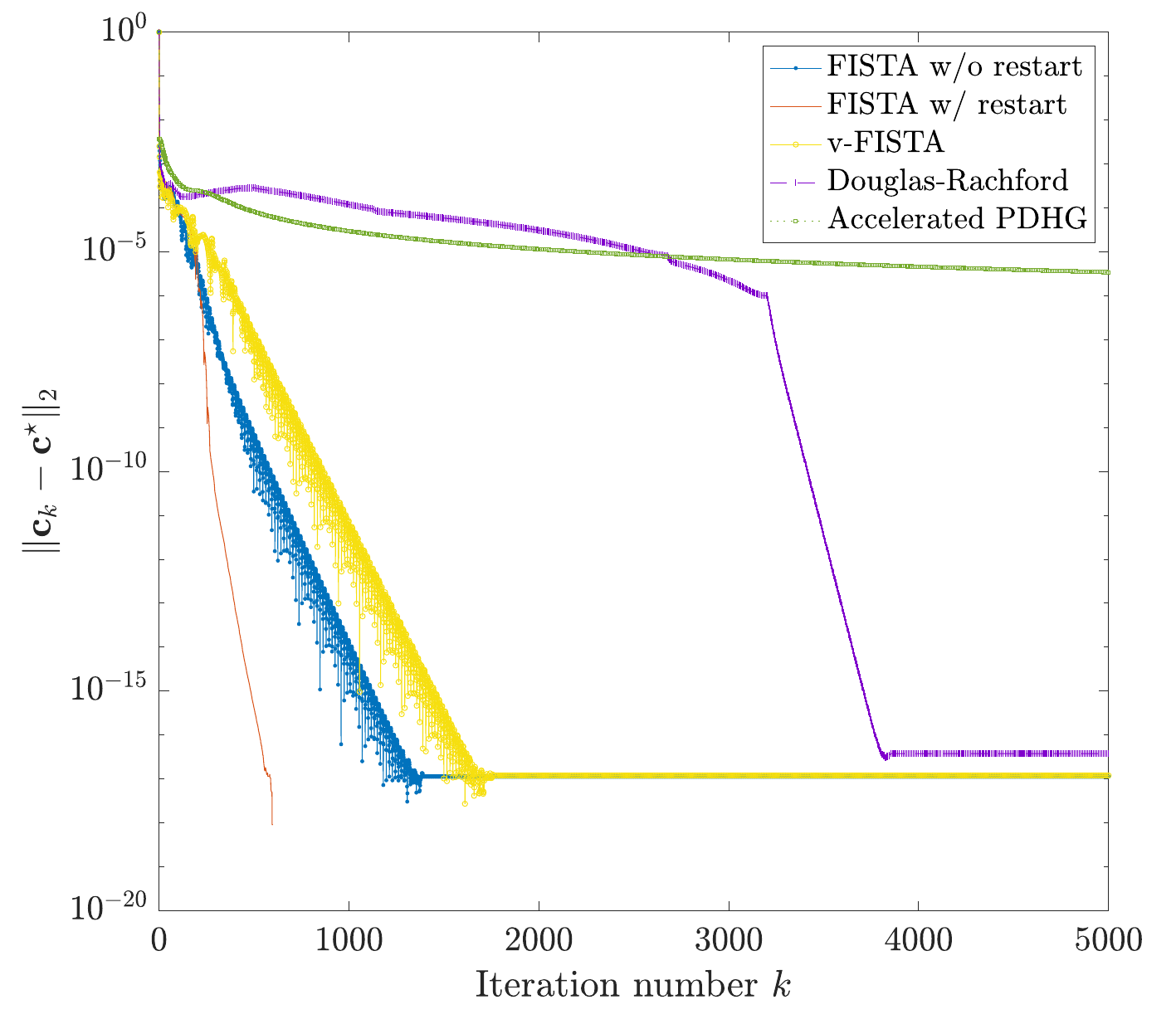}
  \caption{The comparison of convergence curve of several methods for approximating $f_2(x)$ \eqref{equ:tsin} with positivity constraints when $n=20$. The reference minimizer $\mathbf{c}^\star$ is obtained numerically via restarted FISTA with $5,000$ steps. For simplicity, Douglas-Rachford splitting uses the same step size as the FISTA method. We emphasize that Douglas-Rachford splitting could be much faster if tuning parameters.}
\end{figure}

We conduct further examination of the proposed method using this example. As mentioned in Section \ref{sec:method}, FISTA with adaptive restart is not the only choice, so we show the comparison of convergence of several algorithms. The algorithms we compared include accelerated Primal-Dual Hybrid Gradient (PDHG) method \cite{chambolle2011first,pock2011diagonal} for primal-dual problem \eqref{equ:probpd}; projected gradient method, FISTA without restart and Douglas–Rachford splitting \cite{lions1979splitting} for dual problem \eqref{equ:probd}. The schemes used in our experiment are shown in Table \ref{tbl:scheme}. The convergence of these methods for the case $n=20$ is presented in Figure \ref{fig:Tsincompare}.  Asymptotic linear convergence is observed for all the methods with different rates except accelerated PDHG. Our proposed method, FISTA with restart converges to the round-off error in around $600$ iterations with the highest asymptotic rate. The effectiveness of adaptive restart could be observed through comparison to FISTA without restart. \RV{The v-FISTA method behaves even worse than FISTA without restart since the system $\mathbf{B}\K^\dagger\mathbf{B}^T$ is ill-conditioned as mentioned in Remark \ref{rmk:vfista}, with condition number $\kappa \sim \mathcal{O}(10^{17})$.} \RV{We emphasize that in this test the parameter (step size $\eta$) is not tuned for Douglas-Rachford splitting, which could be much faster with tuned parameters.}

We are also interested in the number of iterations needed for convergence, and the influence of the following parameters on the performance:
\begin{itemize}
    \item Number of constraints $C$;
    \item Number of samples used for approximation $K$;
    \item Polynomial order $n$;
    \item Dimension $d$.
\end{itemize}
In this one-dimensional example, we design tests for the first three parameters. The tests for dimensions will be presented in the next subsection and Figure \ref{fig:hddim}.

\begin{figure}[htbp]
  \centering
  \label{fig:Tsin_it}
  \includegraphics[width=.49\textwidth]{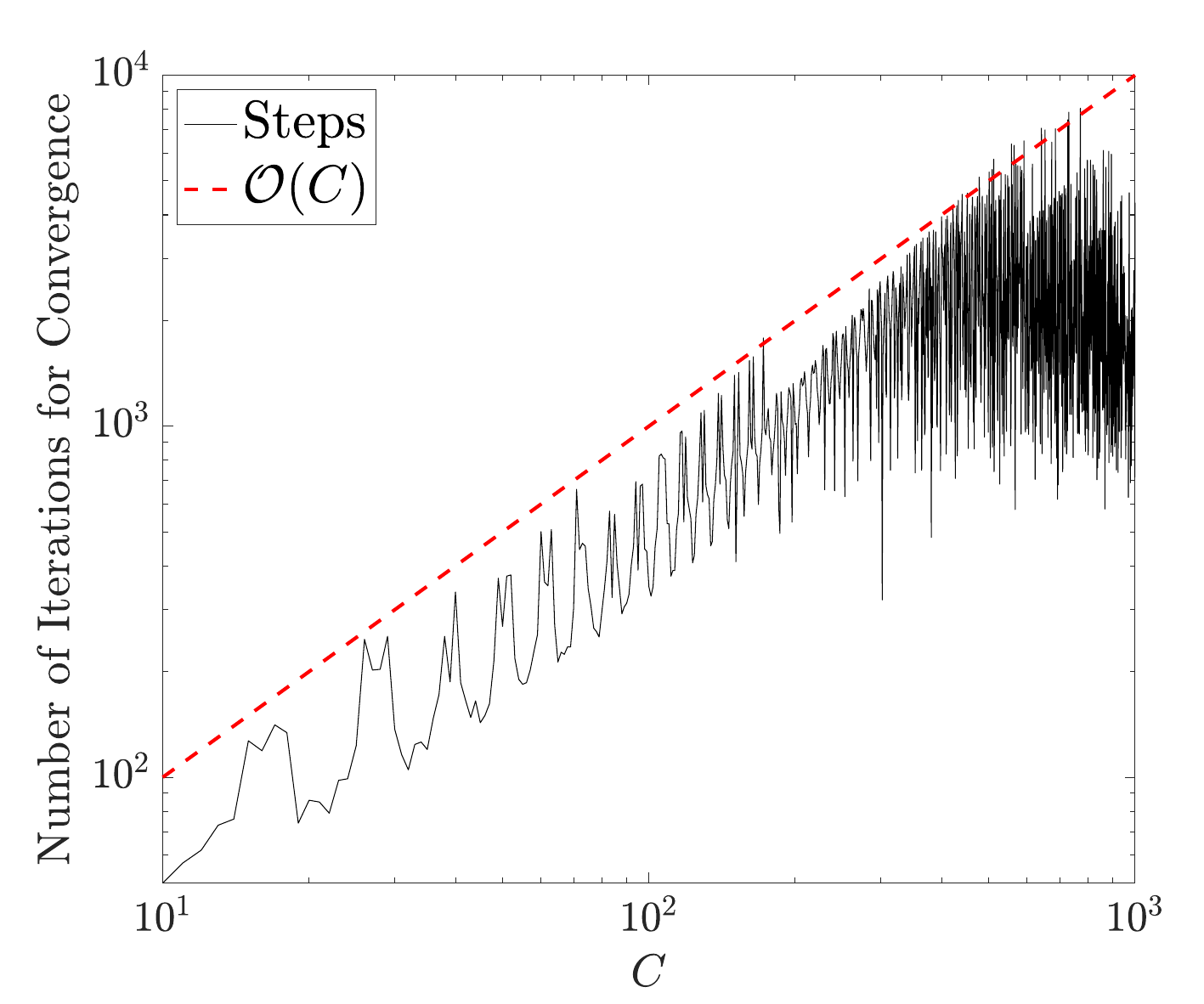}
  \includegraphics[width=.49\textwidth]{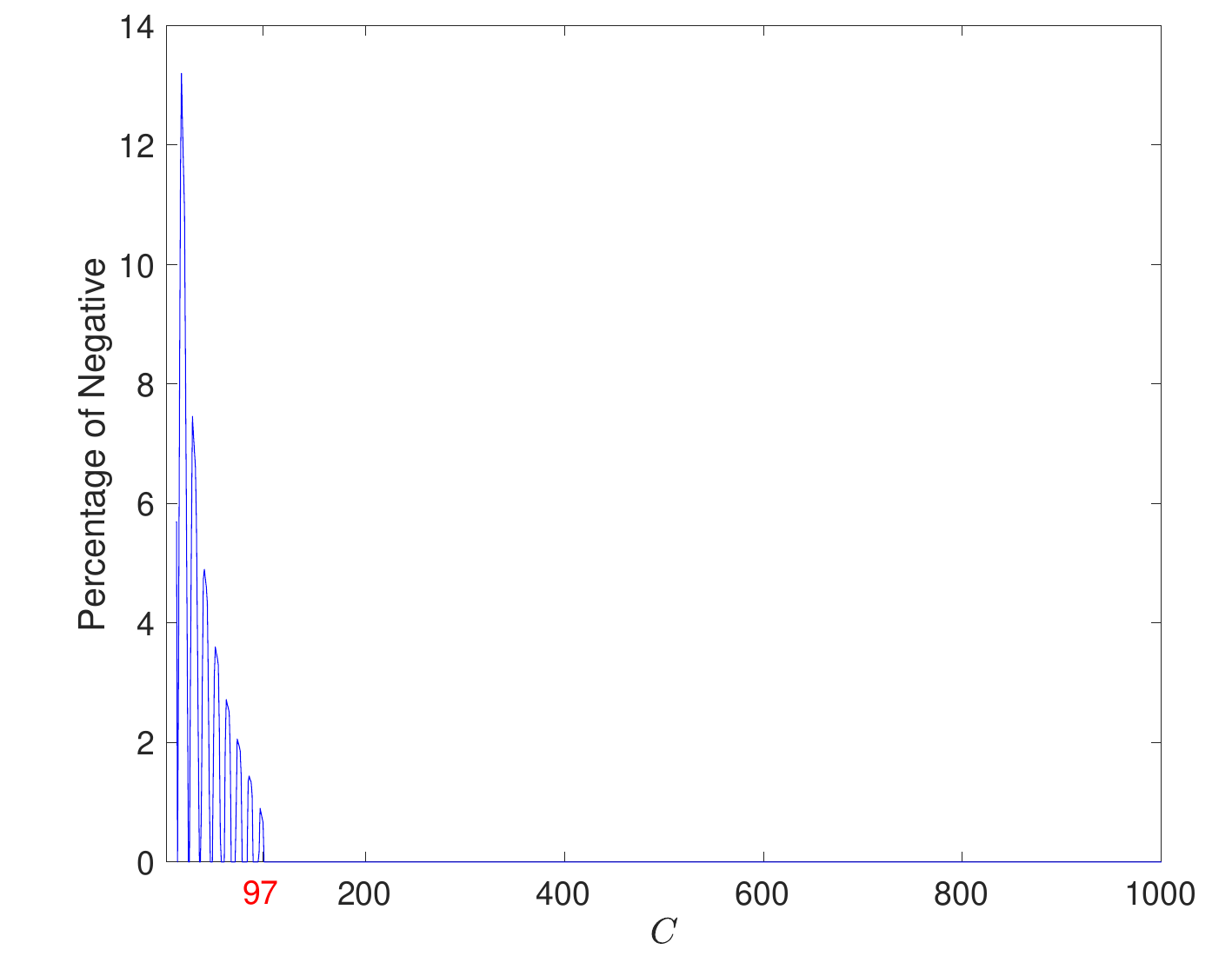}
  \includegraphics[width=.49\textwidth]{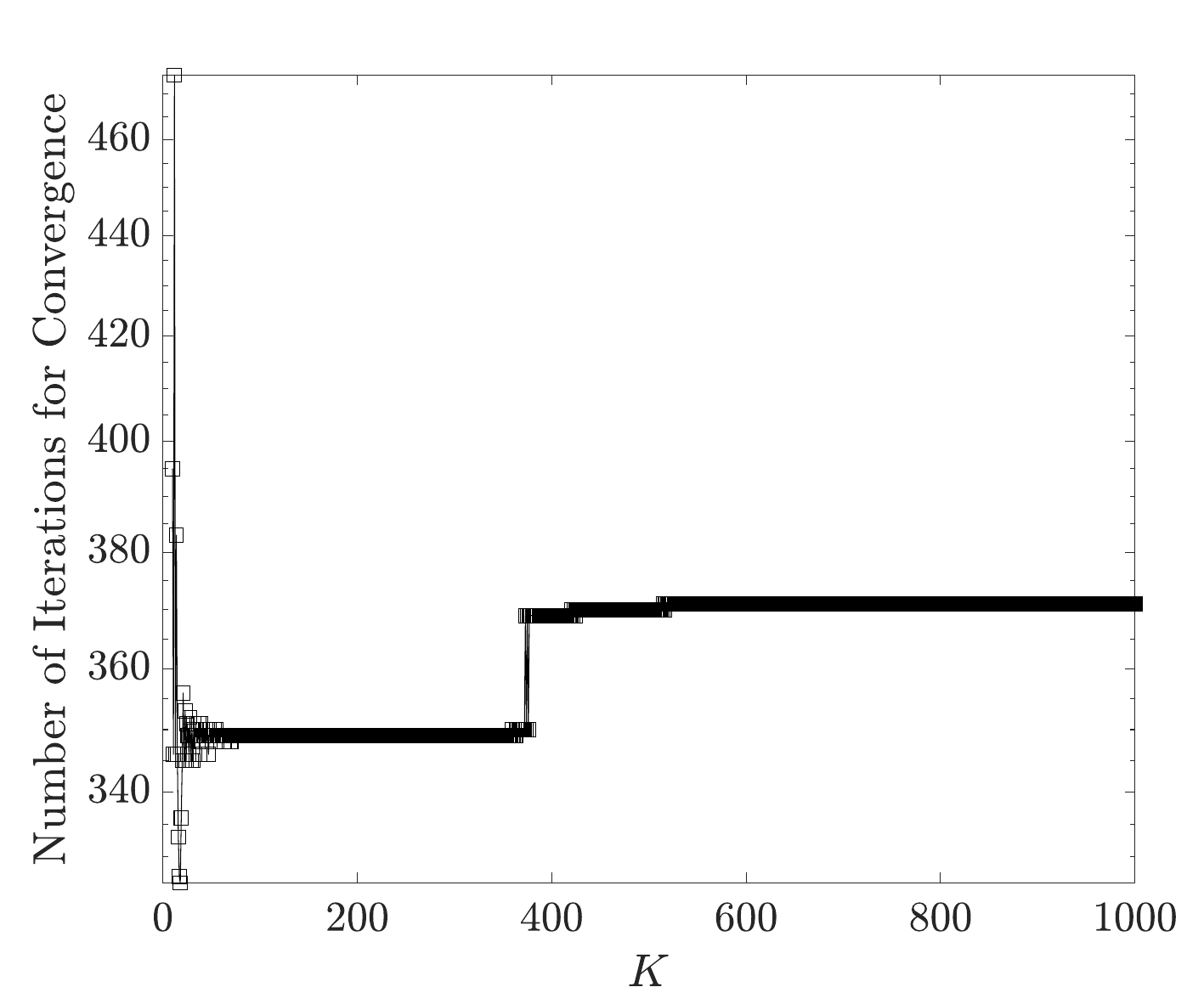}
  \includegraphics[width=.48\textwidth]{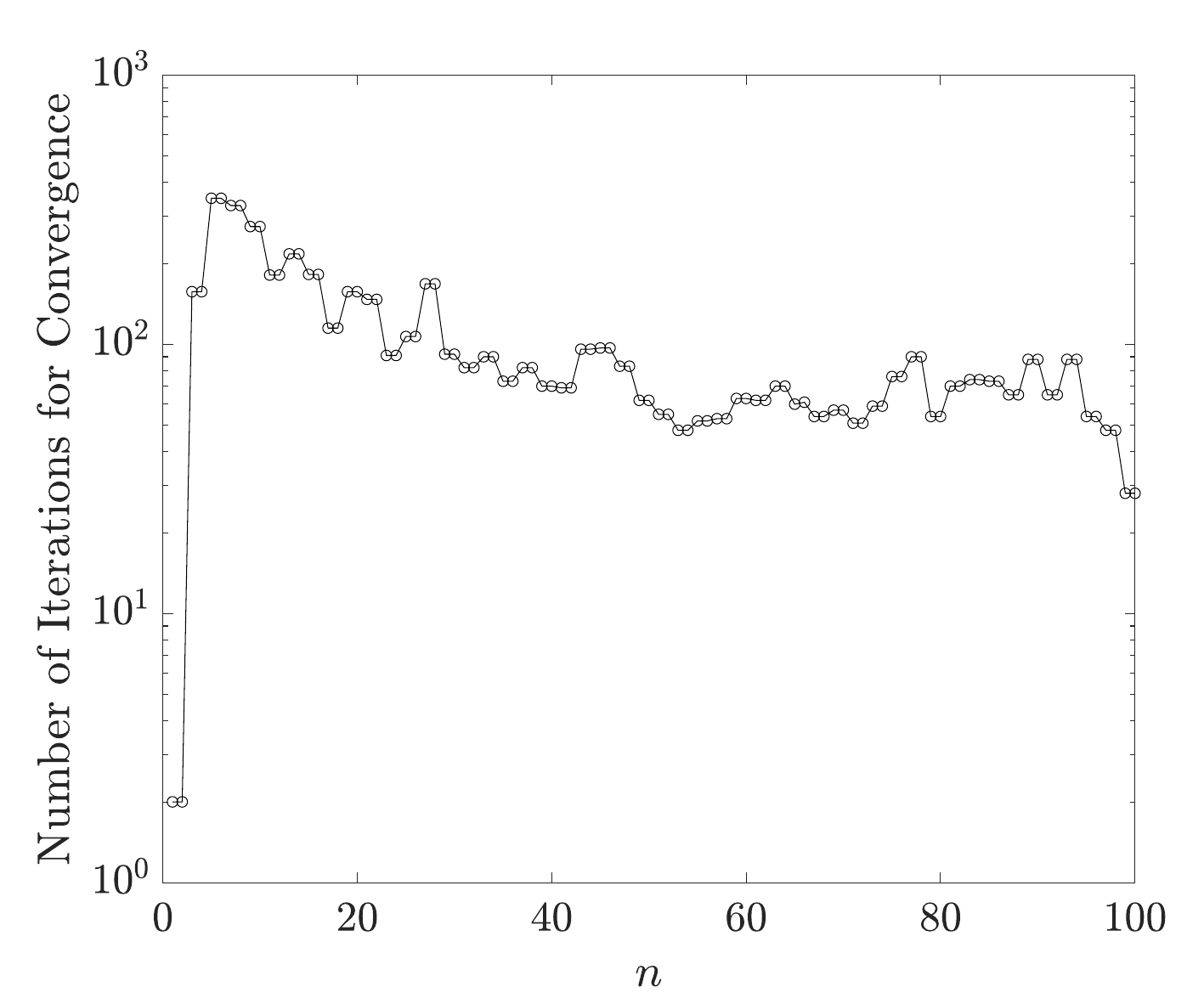}
  \caption{The required iteration number for restarted FISTA to converge, with respect to the number of constraints $C$ (upper left), the number of approximation samples $K$ (lower left), and polynomial order $n$ (lower right) for finding a positive approximation to the function $f_2(x)$ \eqref{equ:tsin}. The percentage of negative sample points in $10^4$ independent randomly tested samples with respect to $C$ is shown in upper right.}
\end{figure}

We record the number of iterations needed for restarted FISTA to converge for each parameter. The results are shown in Figure \ref{fig:Tsin_it}. For $C$, we keep $n=5$, $K=50$ and test different values of $C$ from $10$ to $1000$. As $C$ increases,   numerically the number of iterations needed scales like $\mathcal{O}(C)$. For the other two parameters $K$ and $n$, there is no significant growth of iteration numbers as these parameters are enlarged. These tests show that the parameter $C$ affects the performance the most. \RV{This can be partially explained from a close look of \eqref{equ:probd2}, where $C$ affects the dimension of the optimization problem while $K$ and $n$ only affect the condition.}

Then another question is how would $C$ affect the effectiveness of the positive approximation. In other words, the proposed method only strictly enforces the positivity at $C$ chosen points. As $C$ increases, we expect the minimizer to \eqref{primal} will give a polynomial approximation producing \RV{fewer points where  negative values emerge}.
We check the approximation polynomial values at $10^4$ equidistant points on $[-1,1]$, and compute the percentage of negative point values. For example, assume these test samples are denoted by $\z_i$, $i=1,...,L$, then we compute $|\{\z_i, \widetilde f(\z_i)\geq 0\}|/L$ for  different $C$ from $10$ to $1000$. The result is shown in the upper right of Figure \ref{fig:Tsin_it}. We observed a rapid descent of this percentage. The number of negative points concentrates at $0$ for $C>97$ which indicates the polynomial is positive mostly on the whole domain.

\subsubsection{A step function}
\begin{figure}[htbp]
  \centering
  \label{fig:si}
  \includegraphics[width=.49\textwidth]{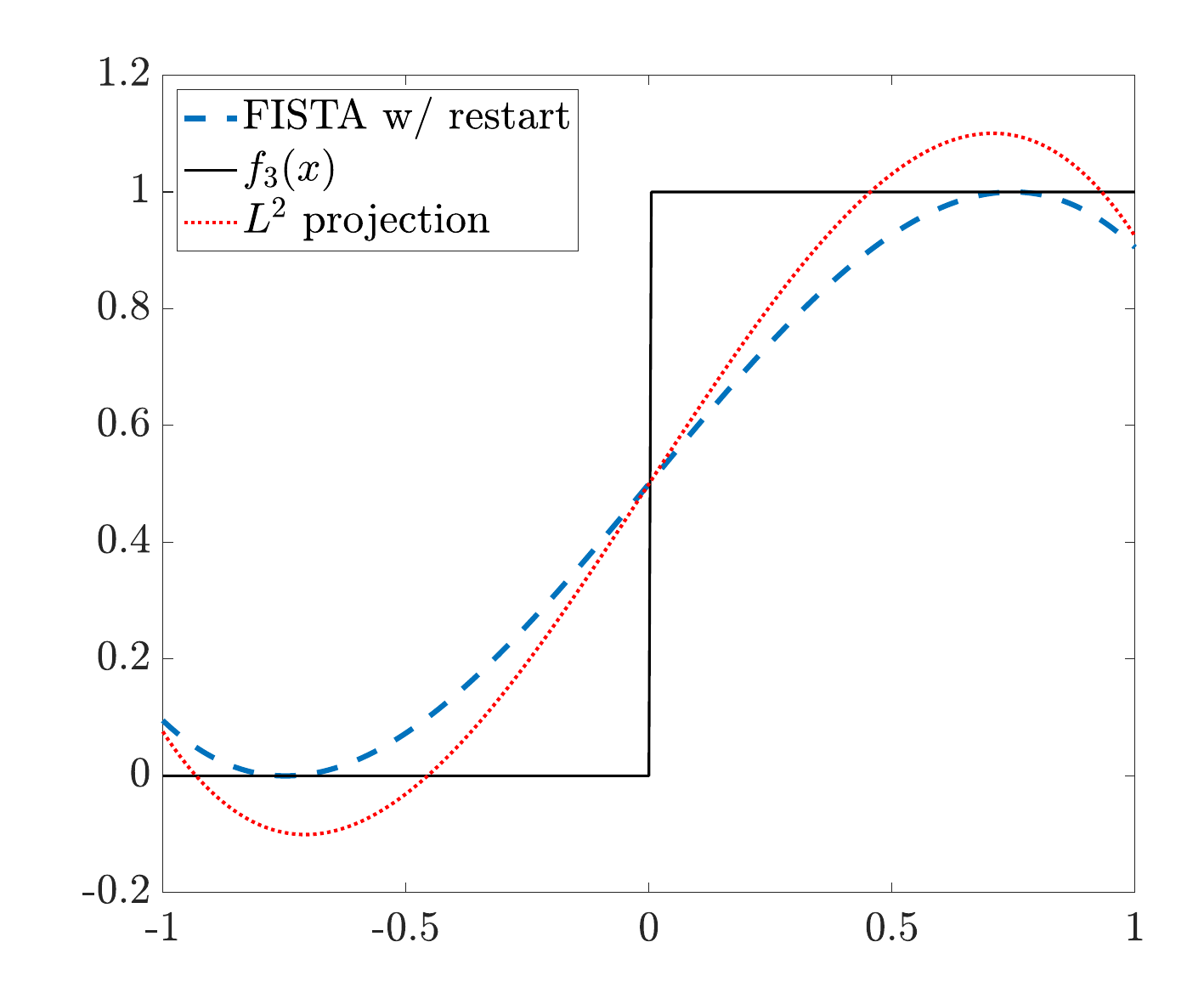}
  \includegraphics[width=.49\textwidth]{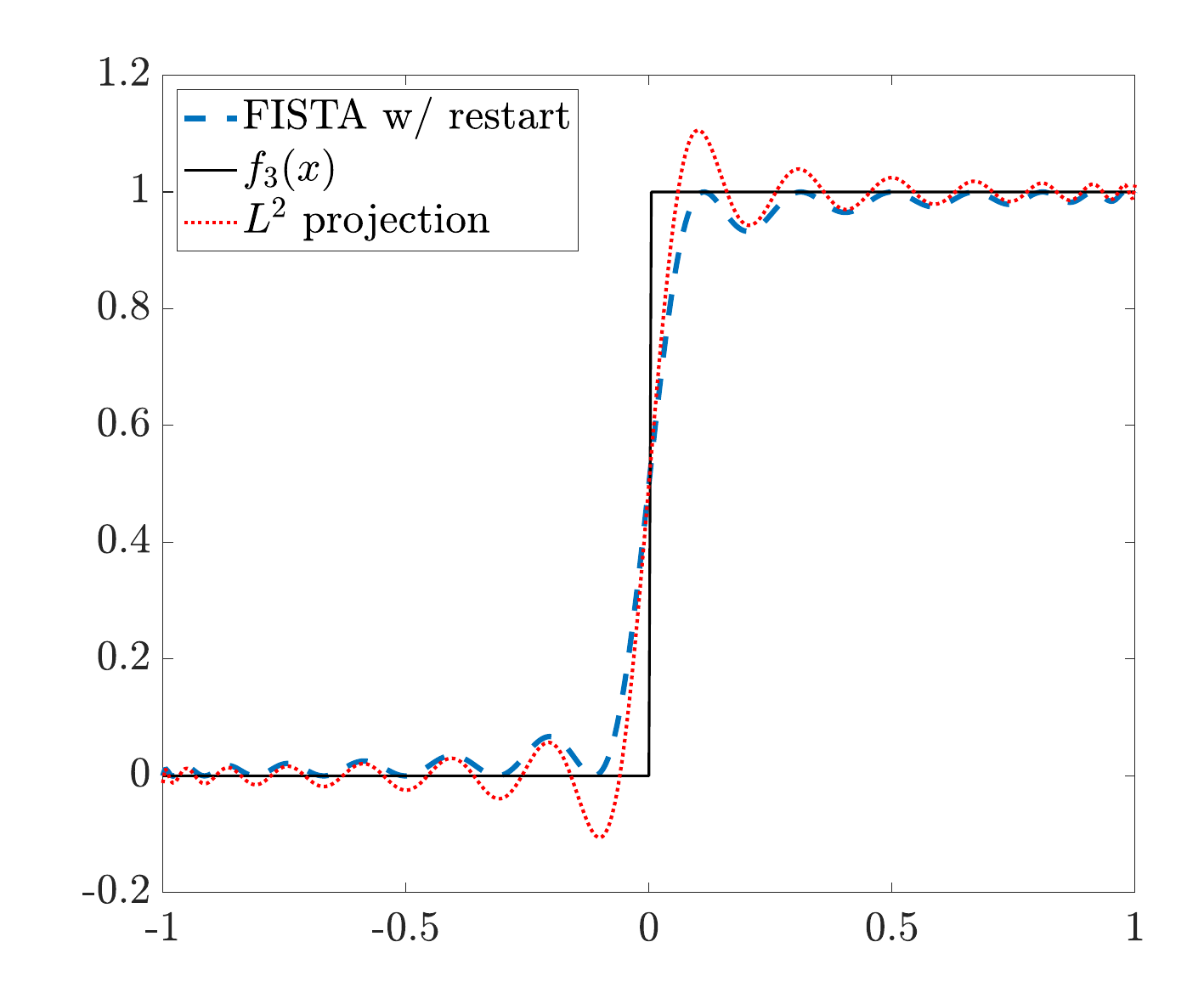}
  \caption{The approximated polynomial for $f_3(x)$ \eqref{equ:f3} with boundedness constraints for $n=5$ (Left), $n=30$ (Right). The bounded approximation is found via solving \eqref{primal} by restarted FISTA on \eqref{equ:probd2}}
\end{figure}

We consider a step function. The polynomial approximation of this function suffers from oscillations which result in overshoots and undershoots. So we consider a bound-preserving approximation using the proposed method. To achieve it, $M=251$ equidistant grid points are used for bound-preserving constraints. We set $\mathbf{B}=\begin{bmatrix}\mathbf{\Psi}_p\\-\mathbf{\Psi}_p\end{bmatrix}$ and $\mathbf{b}=\begin{bmatrix}\mathbf{0}+\epsilon\\\mathbf{-1}+\epsilon\end{bmatrix}$. This enforces the condition $ 0 < \widetilde f(\y_i) < 1$ for each $\y_i$, $i=1,...,M=251$. The results for $n=5$ and $n=30$ are shown in Figure \ref{fig:si}. It is observed, our approximation could almost preserve boundedness.

\subsection{Two-dimensional Examples}
In this section, we consider approximating two-dimensional functions. The domain is set to be $D=[-1,1]\times [-1,1]$. We use $31\times 31$ uniform rectangular grids for the approximation points $\x_i$, thus $K=961$. 
\subsubsection{Gaussian peak function}
\begin{figure}[htbp]
  \centering
  \label{fig:Gaussian2d}
  \includegraphics[width=.3\textwidth]{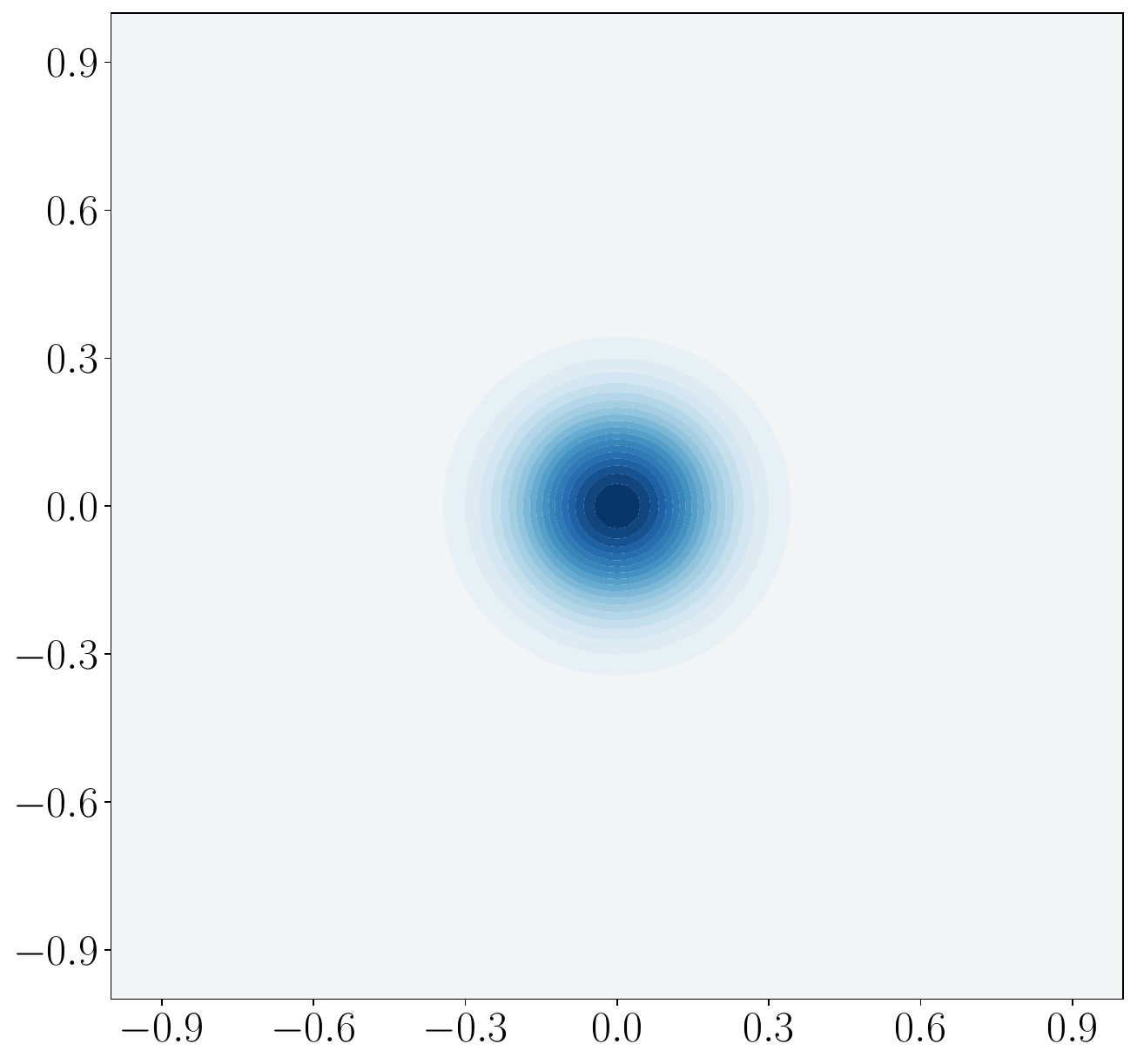}
  \includegraphics[width=.3\textwidth]{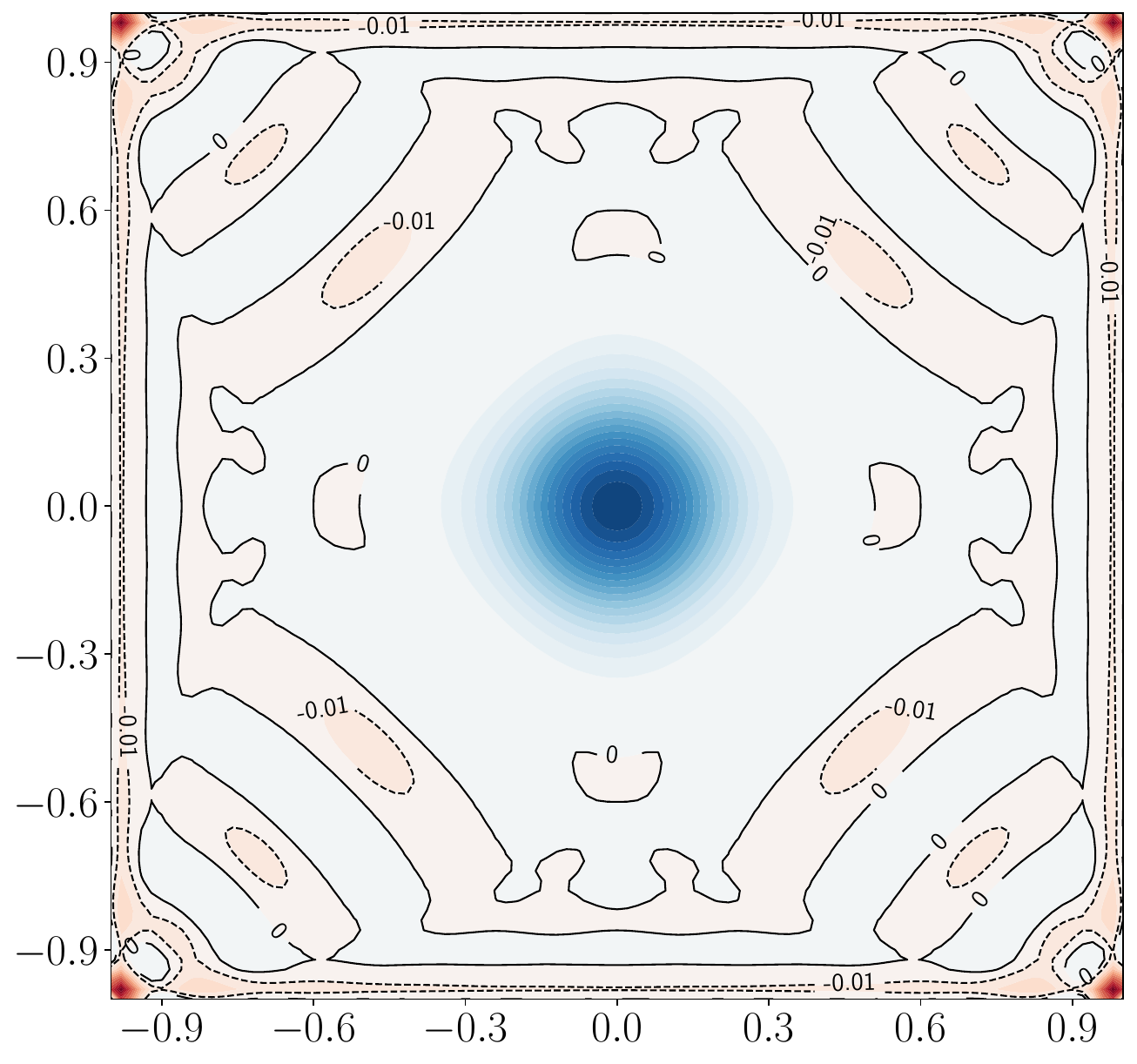}
  \includegraphics[width=.362\textwidth]{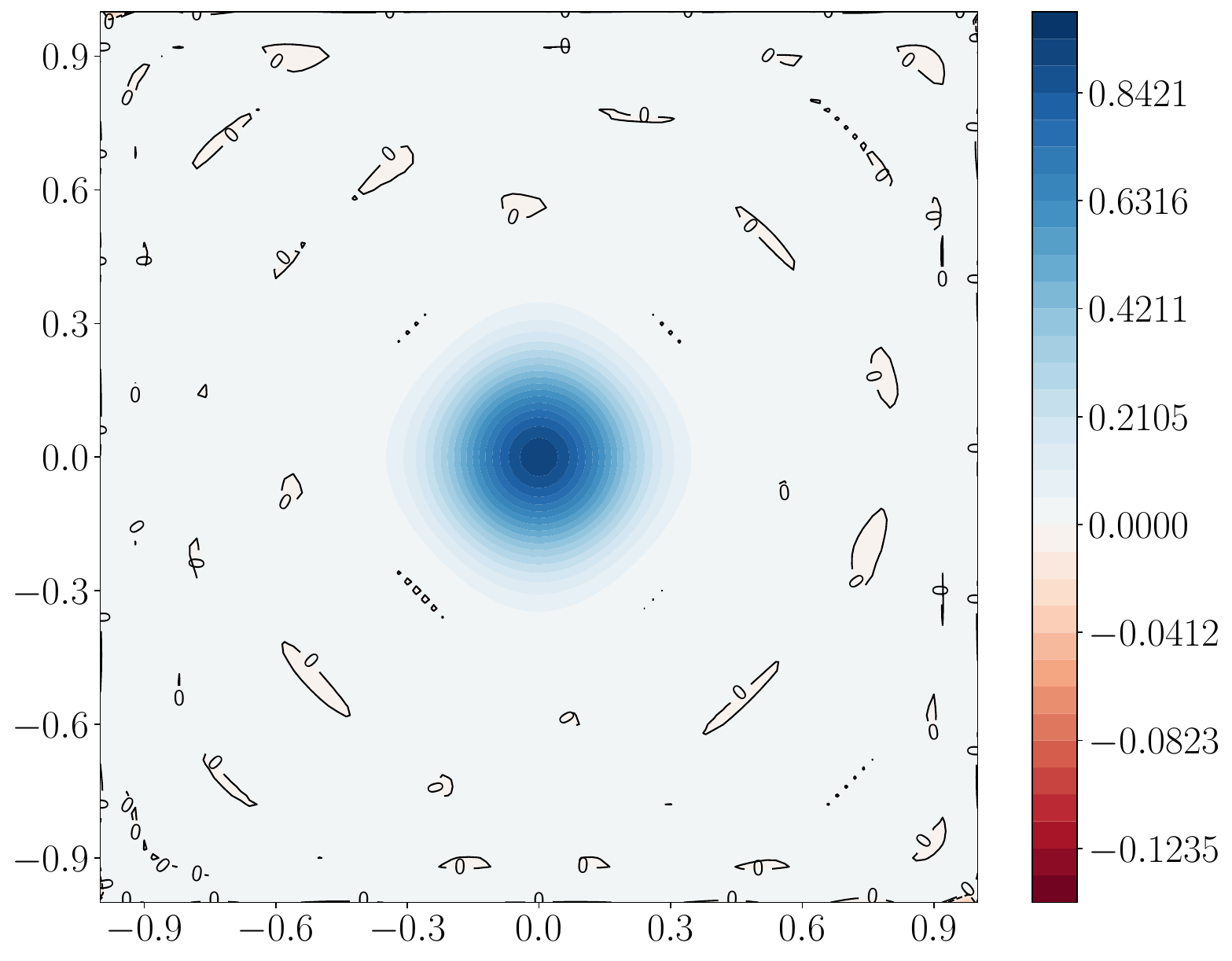}
  \caption{The contour plots of $f_4(\x)$ \eqref{equ:gaussian} (left), $L^2$ projection when polynomial order $n=20$ (Middle) and positive approximation via solving \eqref{primal} by restarted FISTA on \eqref{equ:probd2} (Right). In these graphs, positive values are marked in blue color, negative values are marked in red. A few contour lines are drawn using black solid (zero value) and dashed (a few negative values) lines.}
\end{figure}

We set $\sigma_1=\sigma_2=10$ and $\omega_1=\omega_2=0.5$ in the formula \eqref{equ:gaussian}. A smooth peak could be observed in the contour plot, see Figure \ref{fig:Gaussian2d}. Standard $L^2$ Approximations around this peak will lead to negative oscillations. To construct positivity preserving polynomial approximation, we randomly sample $3,000$ points from uniform distribution and enforce the positive conditions on these points. The results with $n=20$ are shown in Figure \ref{fig:Gaussian2d}. In these contour plots, the color blue and red are used to indicate positive and negative values, respectively. A few contour lines are also shown in the graph using black color. It is clear to see that negative values occur widely around the peak of standard $L^2$ approximations. The approximation of our method shows a significant improvement.

\subsubsection{Continuous peak function}
\begin{figure}[htbp]
  \centering
  \label{fig:Continuous2d}
  \includegraphics[width=.3\textwidth]{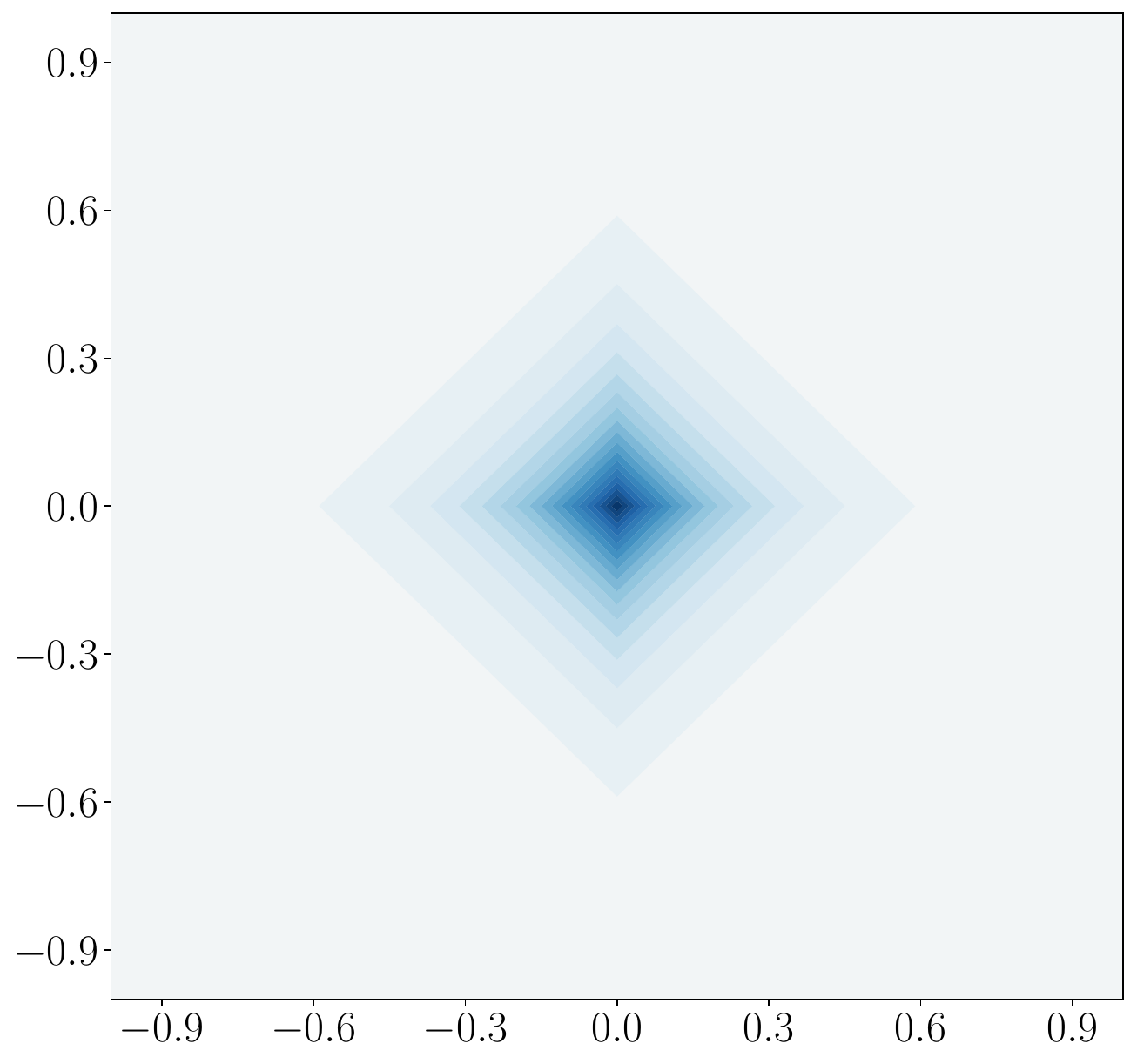}
  \includegraphics[width=.3\textwidth]{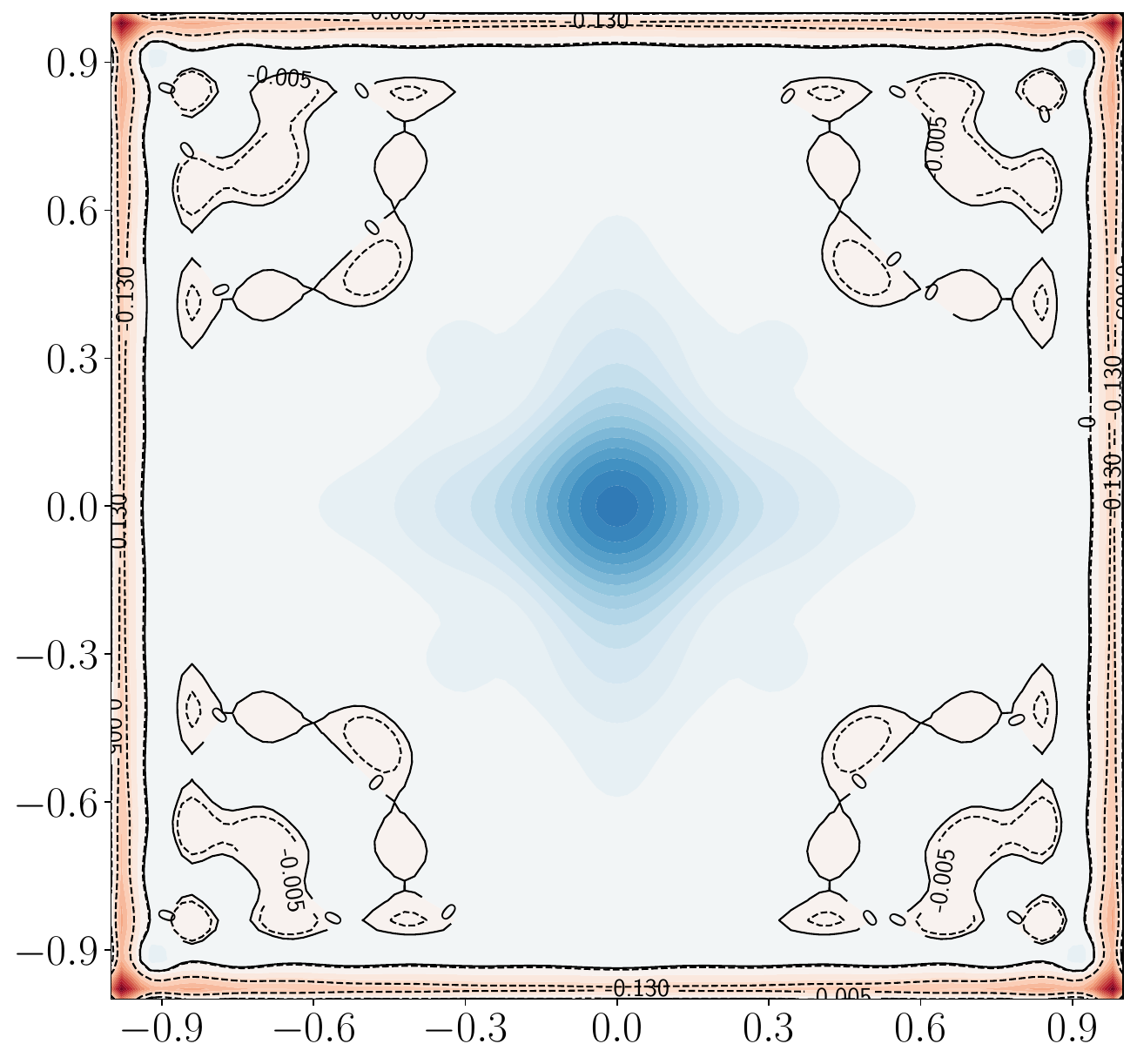}
  \includegraphics[width=.356\textwidth]{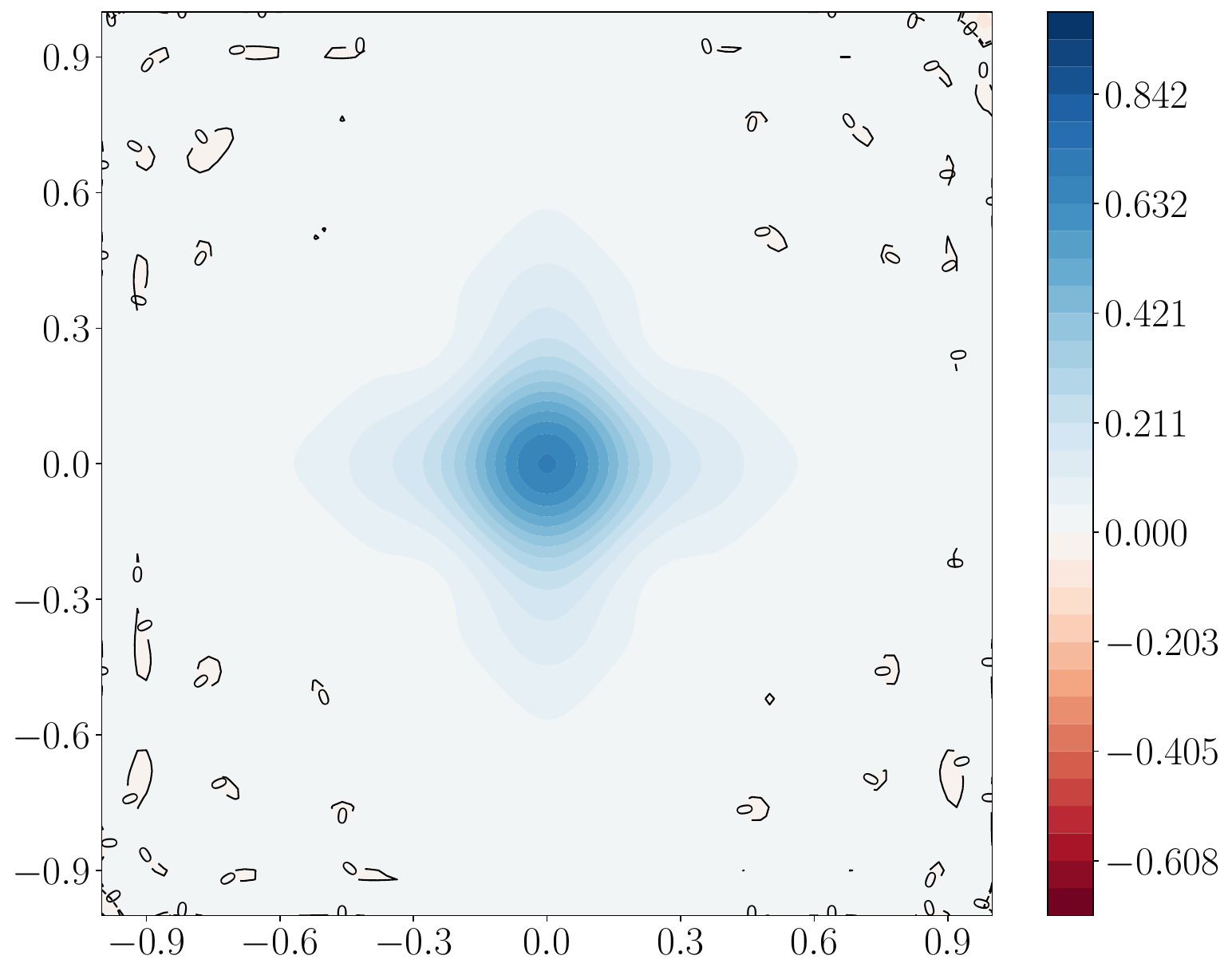}
  \caption{The contour plot of $f_5(\x)$ \eqref{equ:f5} (left), $L^2$ projection when polynomial order $n=20$ (Middle) and positive approximation via solving \eqref{primal} by restarted FISTA on \eqref{equ:probd2} (Right). In these graphs, positive values are marked in blue color, negative values are marked in red. A few contour lines are drawn using black solid (zero value) and dashed (a few negative values) lines.}
\end{figure}

With the same condition as $f_4(\x)$, we tested the performance of $f_5(\x)$. The parameters are similarly set to be $\sigma_1=\sigma_2=10$ and $\omega_1=\omega_2=0.5$. This function has a non-smooth peak around the origin. The contour plots of the function $f_5(\x)$, $L^2$ approximation, and the approximation using our methods are shown in Figure \ref{fig:Continuous2d}. We observed negative values for standard $L^2$ projection around the four edges of the domain, especially on the four corners. This is significantly improved with our method.

\subsubsection{Corner peak function}
\begin{figure}[htbp]
  \centering
  \label{fig:Cornerpeak2d}
  \includegraphics[width=.3\textwidth]{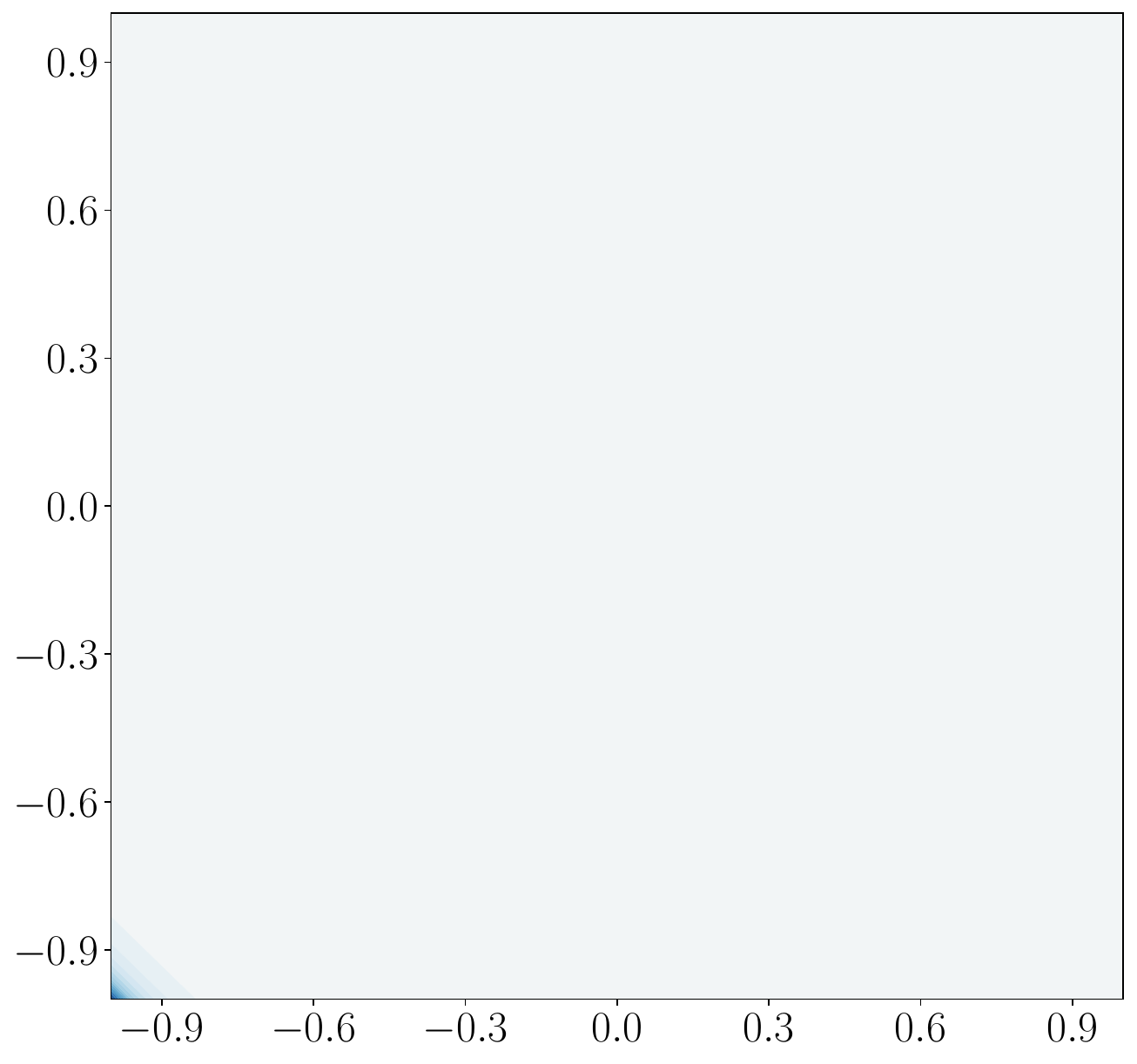}
  \includegraphics[width=.3\textwidth]{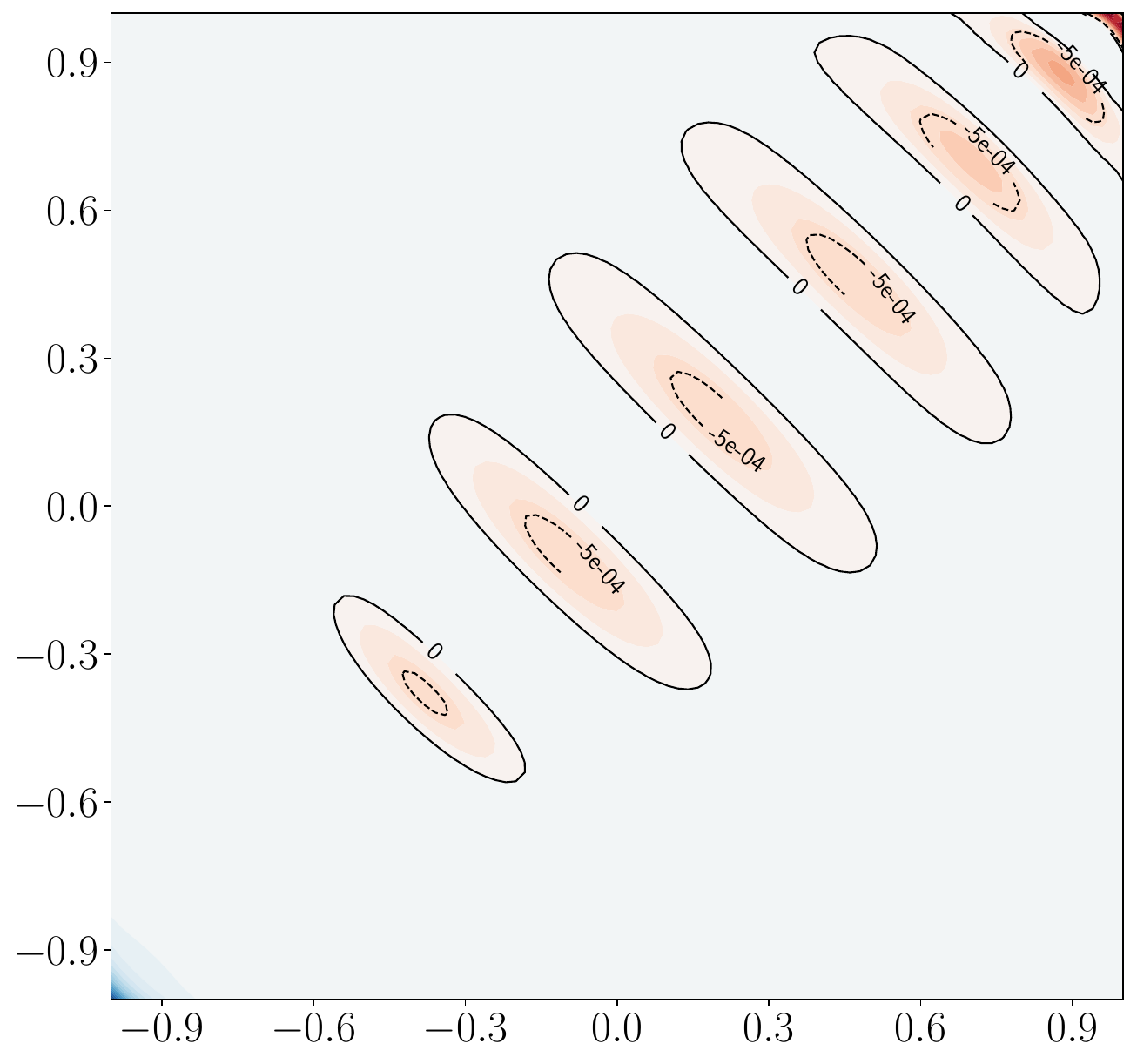}
  \includegraphics[width=.362\textwidth]{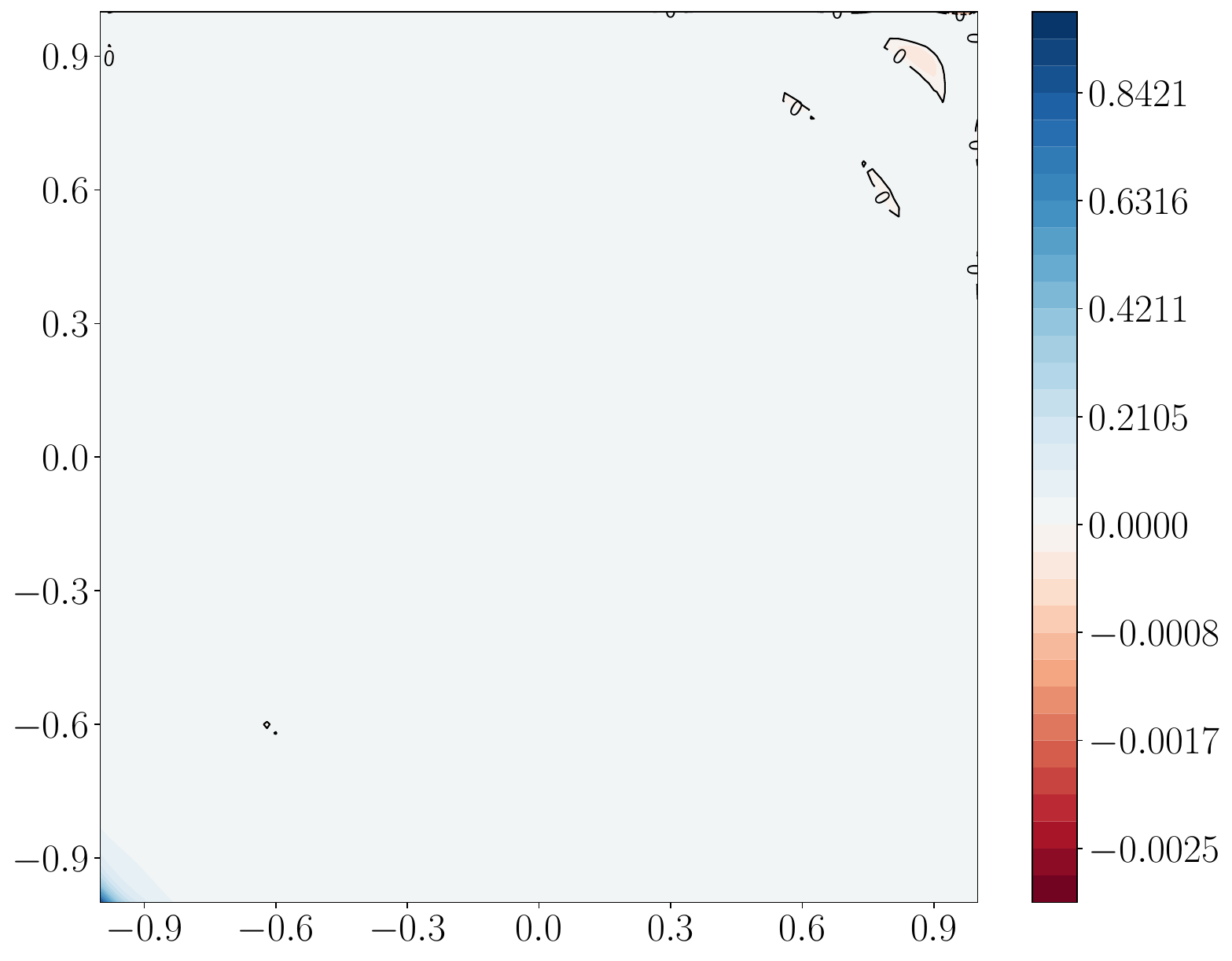}
  \caption{The contour plots of $f_6(\x)$ \eqref{equ:f6} (left), $L^2$ projection when polynomial order $n=20$ (Middle) and positive approximation via solving \eqref{primal} by restarted FISTA on \eqref{equ:probd2} (Right). In these graphs, positive values are marked in blue color, negative values are marked in red. A few contour lines are drawn using black solid (zero value) and dashed (a few negative values) lines.}
\end{figure}

We now consider the example $f_6(\x)$ with parameters $\sigma_1=\sigma_2=20$. This function is called "corner peak" because its value is concentrated in the lower-left corner of the domain, see the left of Figure \ref{fig:Cornerpeak2d}. With standard $L^2$ approximation, it exhibits a large area of negative values along the diagonal of the domain $D$ (Middle of Figure \ref{fig:Cornerpeak2d}). This can be improved significantly by the proposed method, with positivity enforcement on a set of $3,000$ randomly sampled points from $\mathcal{U}(D)$. The result is on the right of Figure \ref{fig:Cornerpeak2d}. 

\subsection{High dimensional Examples}
We present results in dimensions $d \geq 3$, including the immediate dimension $d=10$ and high dimensions $d=100$, $d=200$. We only focus on the Gaussian function $f_4(\x)$ on domain $D=[-1,1]^d$.  We will examine the results with several metrics. One is the approximation error in $L^2$ norm. That is,
\be \label{equ:approxL2}
   \|f-\widetilde f\|_{L^2(D)} \approx \sqrt{\frac{1}{L}\sum_{i=1}^{L} (f(\z_i)-\widetilde f(\z_i))^2},
\ee
where $\z_i$, $i=1,\dots,L$, are sampled uniformly in the domain. These points are different from the ones used for function approximation.
Another metric is the percentage of negative points, which has been defined and used in one-dimensional examples. 

\subsubsection{Dimension 10}
We consider the case  $d=10$. The parameters of $f_4(\x)$ are set as $\sigma_i=10$, $i=1,2, \dots, 10$. To approximate this function, we use $K=2,000$ random generated samples points.
The orders of the polynomial approximations are $n=3,4,5,6,7,8$, with the corresponding cardinality of the polynomial spaces $N=286$, $1,001$, $3,003$, $8,008$, $19,448$ and $43,758$, respectively. The points to enforce positivity are randomly sampled with sizes
from $C=20$ to $12,000$.
To evaluate the performance, we sample another $L=5,000$ points to compute the approximated $L^2$ error and percentage of negative points. 
\begin{figure}[htbp]
  \centering
  \label{fig:d10apppos}
  \includegraphics[width=.49\textwidth]{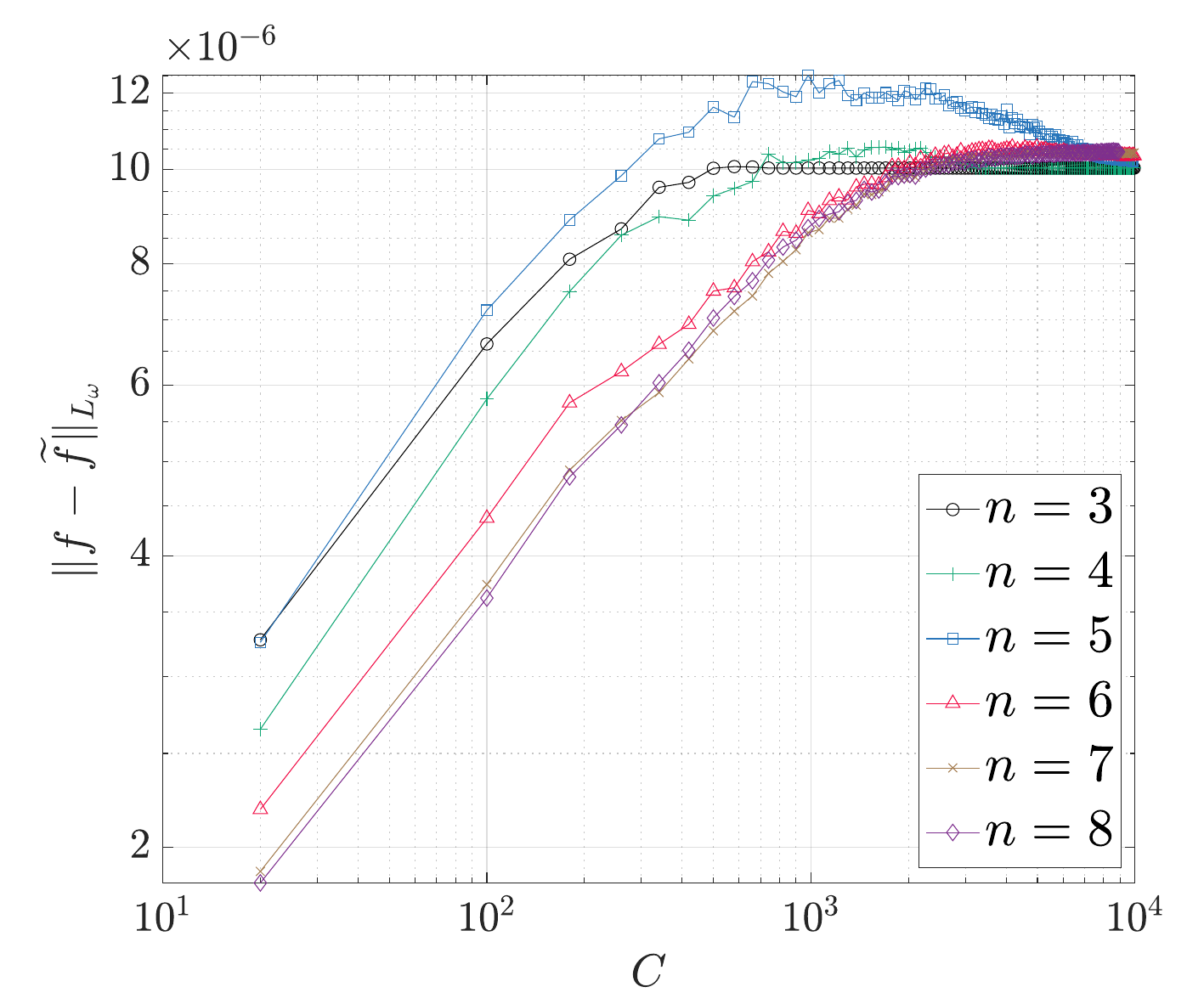}
  \includegraphics[width=.49\textwidth]{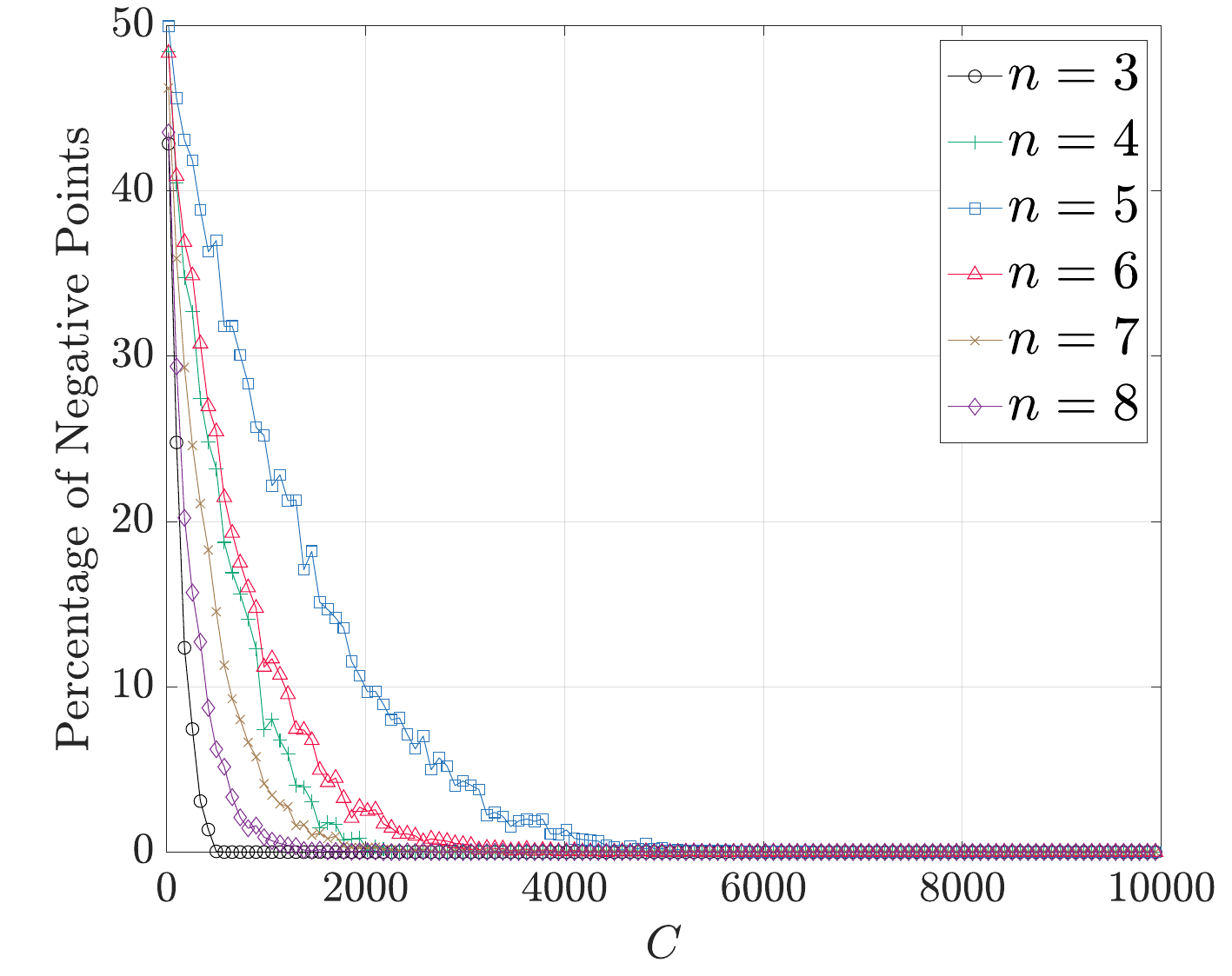}
  \caption{Comparasion of the estimated approximation error \eqref{equ:approxL2} (left) and percentage of negative sample points (right) with respect to the number of constraints $C$ for $10$-dimensional $f_4(\x)$ \eqref{equ:gaussian} using restarted FISTA at polynomial orders $n=3,4,5,6,7,8$. These quantities are computed based on $5,000$ independent randomly tested samples.}
\end{figure}
The results are shown in Figure \ref{fig:d10apppos}. It can be seen that the percentage of the negative points converges to $0$ as the number of positivity enforcing points $C$ increases. The approximation errors do not vary too much with respect to $C$, staying within the range of $2\times 10^{-6} \sim 1.2\times 10^{-5}$. Its magnitude first grows as $C$ increases, and then converges as $C$ becomes larger. For large $C$, we do not observe a significant difference among the approximation errors at the chosen polynomial orders.

\begin{figure}[htbp]
  \centering
  \label{fig:d10conv}
  \includegraphics[width=.99\textwidth]{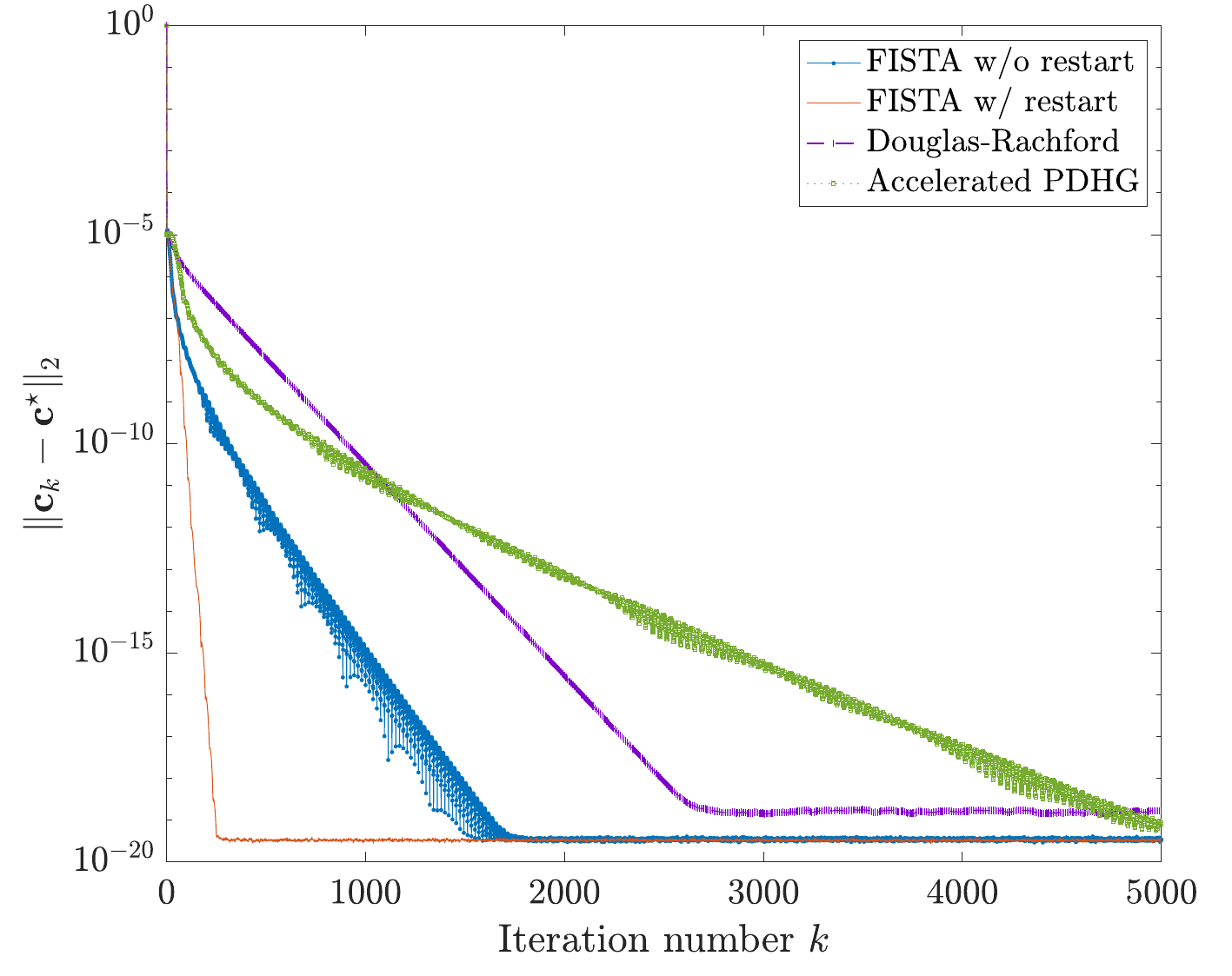}
  \caption{The comparison of convergence curve of several methods for approximating $10$ dimensional $f_4(\x)$ \eqref{equ:gaussian} with positivity constraints when $n=3$. The reference minimizer $\mathbf{c}^\star$ is obtained numerically via restarted FISTA with $5000$ steps. For simplicity, Douglas-Rachford splitting uses the same step size as the FISTA method. We emphasize that Douglas-Rachford splitting could be much faster if tuning parameters. }
\end{figure}

\begin{figure}[htbp]
  \centering
  \label{fig:d10k}
  \includegraphics[width=.99\textwidth]{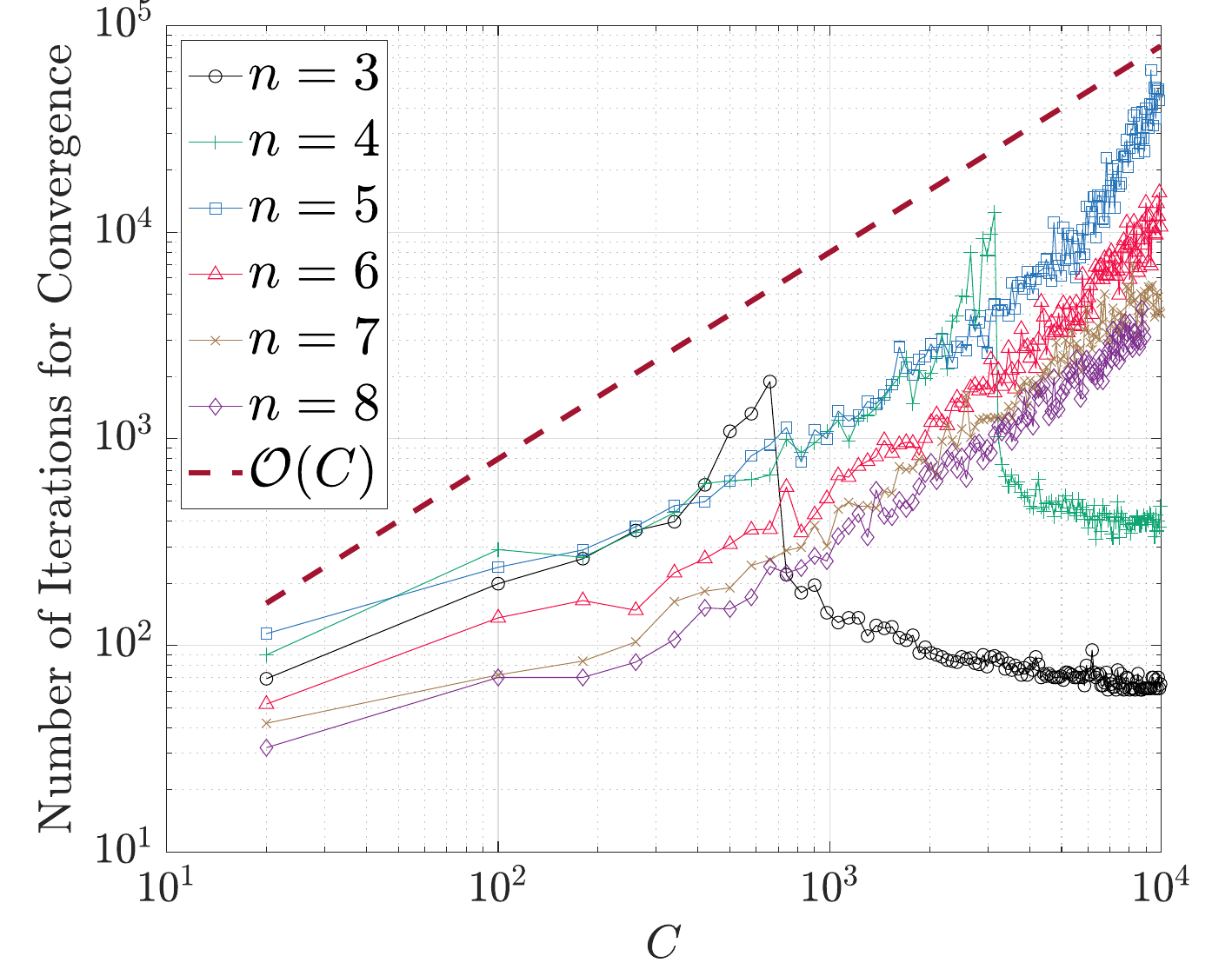}
  \caption{The number of iterations needed for restarted FISTA to converge for $10$-dimensional $f_4(\x)$ \eqref{equ:gaussian} with respect to number of constraints $C$ given polynomial order $n=3,4,5,6,7,8$.}
\end{figure}

In Figure \ref{fig:d10conv}, we show comparison of the proposed method with other splitting algorithms, for the case of $n=3$ and $C=1,000$. We observe superior convergence of the proposed method.
We are also interested in the number of iterations needed for the proposed FISTA with restart, shown in Figure \ref{fig:d10k}. An approximately linear growth w.r.t. $C$ is observed, similar to the one-dimensional case.   It is observed that polynomials with larger orders usually need fewer steps to converge because the condition of matrix $\mathbf{B}$ is better.

\subsubsection{Dimensions of 100 and 200}

\begin{figure}[htbp]
  \centering
  \label{fig:hdapp}
  \includegraphics[width=.49\textwidth]{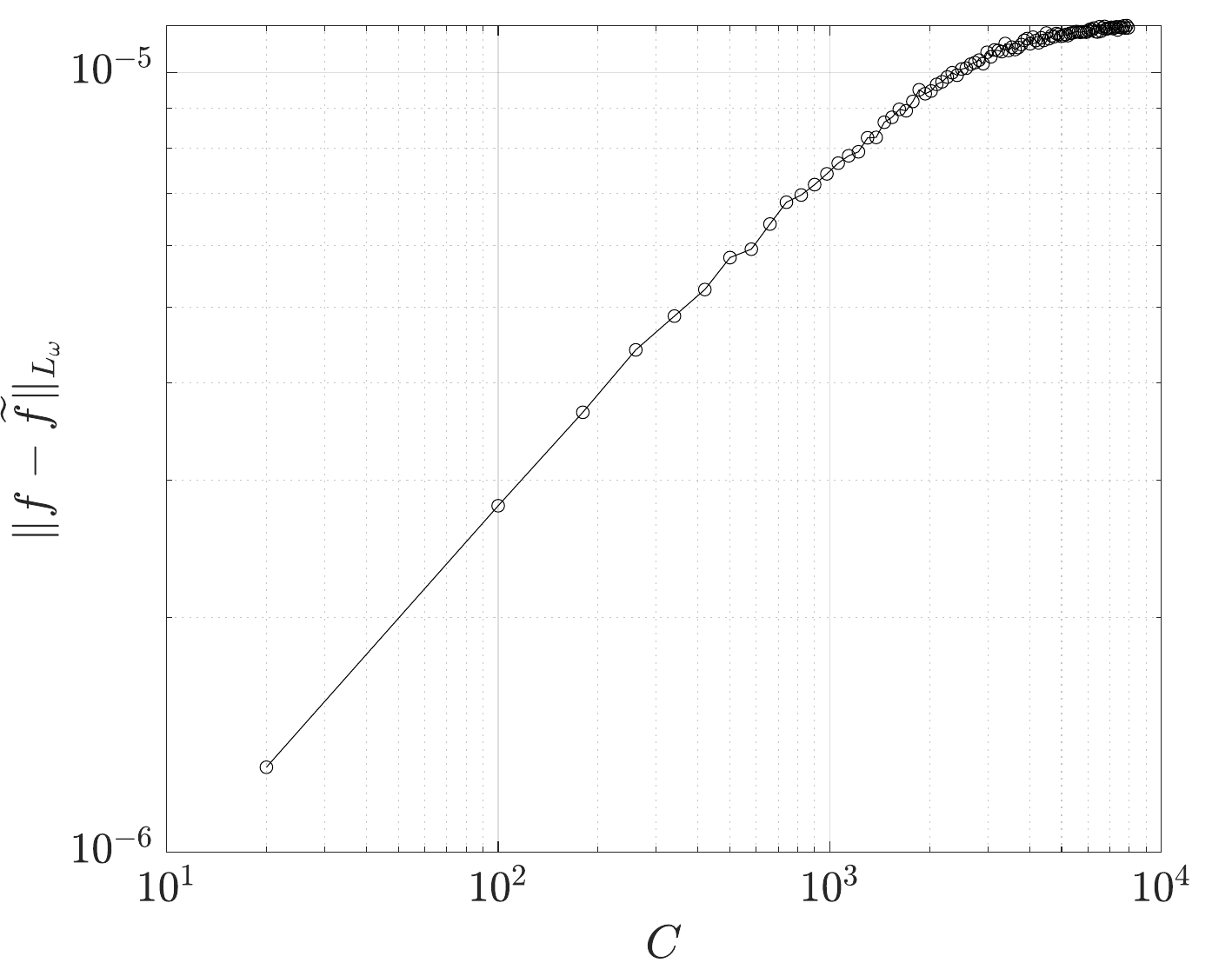}
  \includegraphics[width=.49\textwidth]{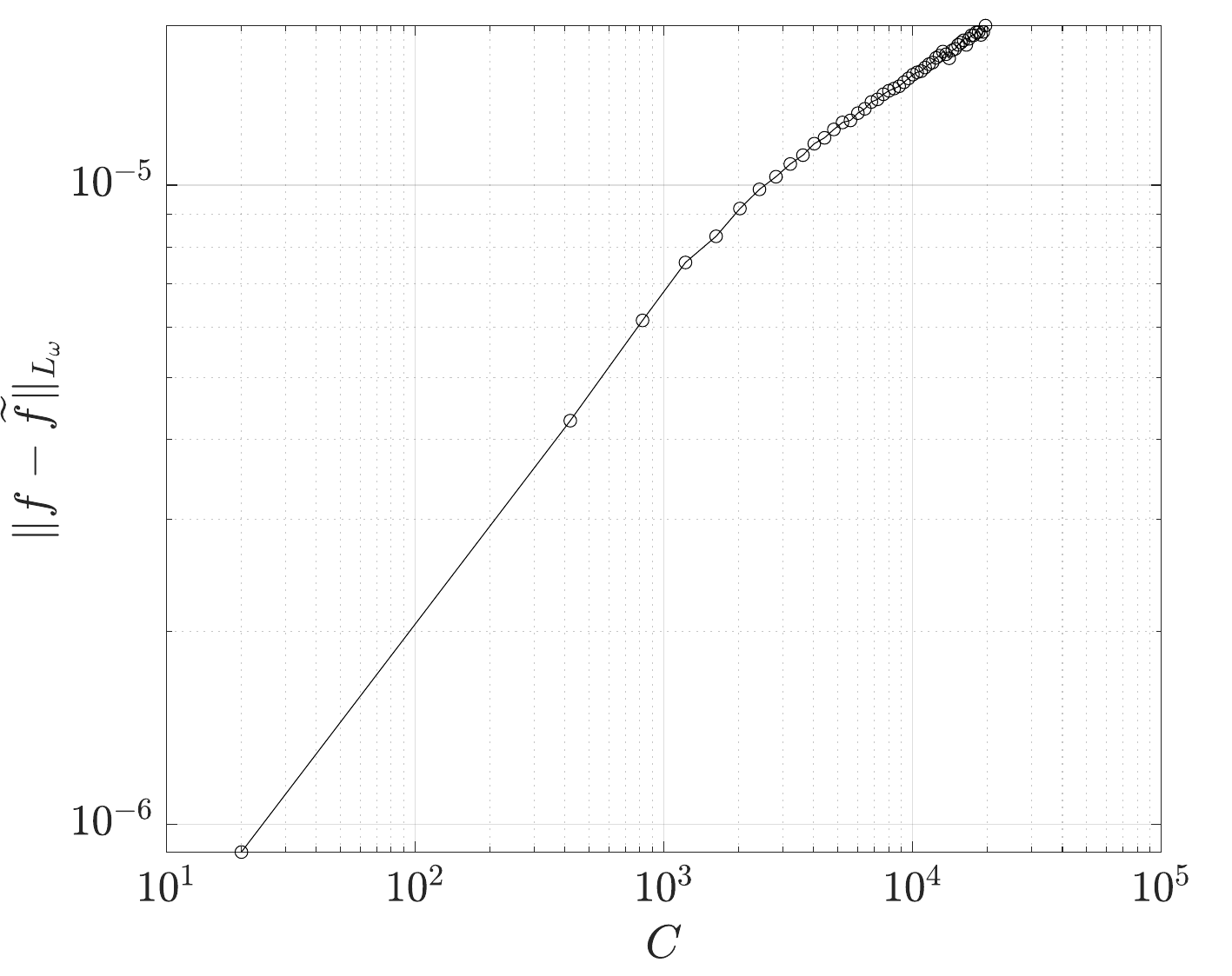}
  \caption{The estimated approximation error \eqref{equ:approxL2} for $100$-dimensional (left) and $200$-dimensional (right) $f_4(\x)$ \eqref{equ:gaussian} with respect to the number of constraints $C$ using restarted FISTA when polynomial orders $n=2$. The errors are computed based on $5,000$ independent randomly tested samples.}
\end{figure}

\begin{figure}[htbp]
  \centering
  \label{fig:hdpos}
  \includegraphics[width=.49\textwidth]{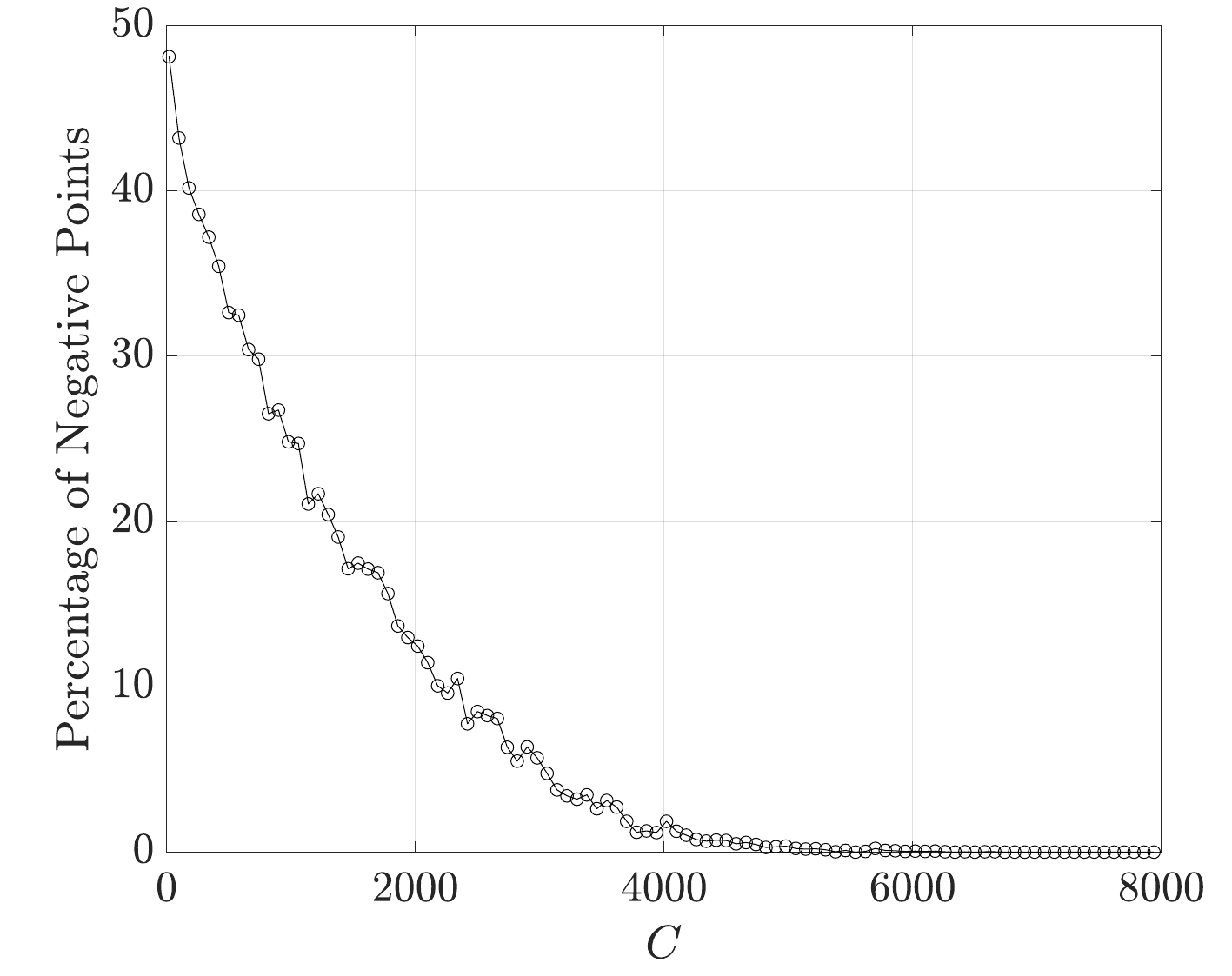}
  \includegraphics[width=.49\textwidth]{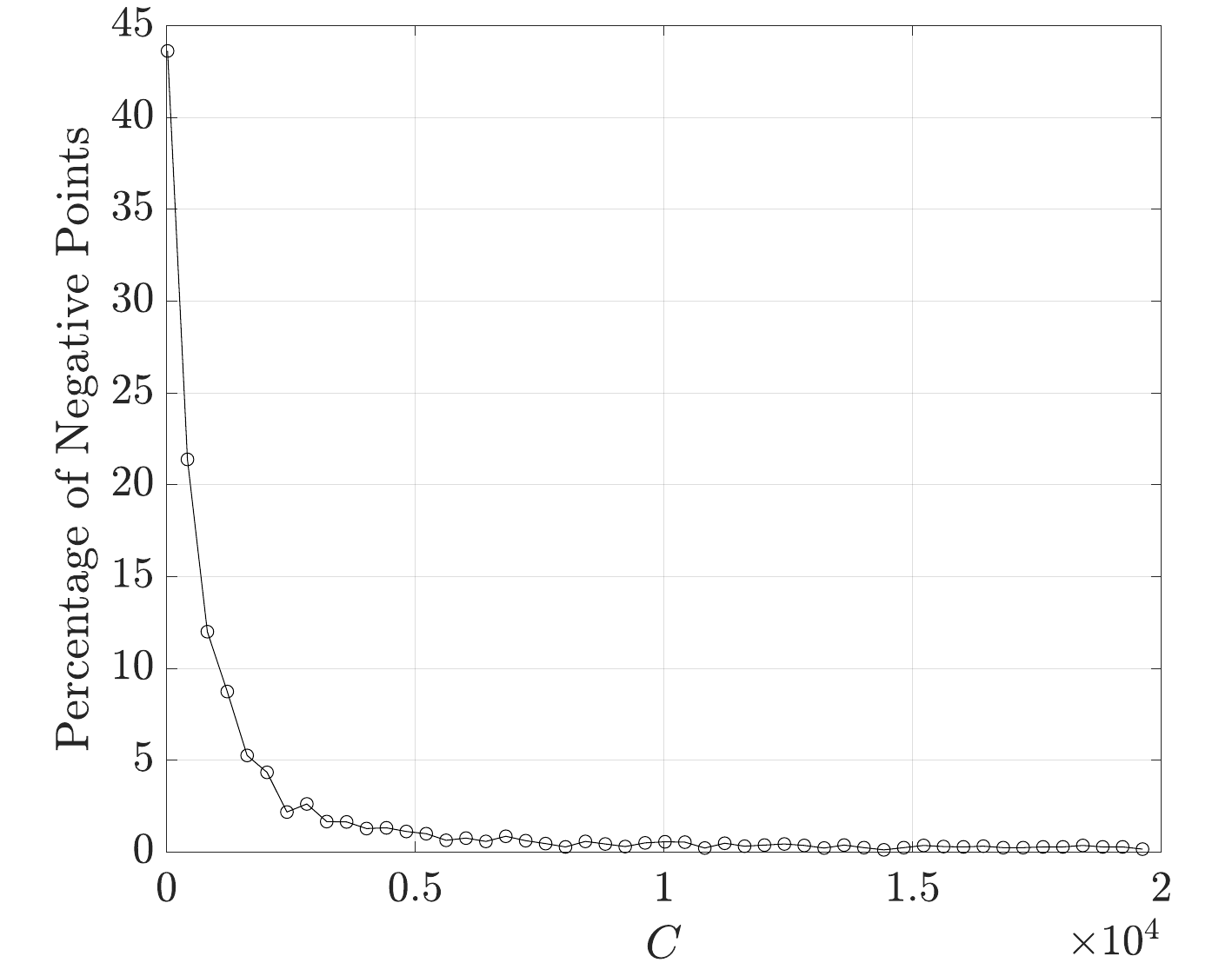}
  \caption{The percentage of negative sample test points of restarted FISTA method (based on $5,000$ independent randomly tested samples) for $100$-dimensional (left) and $200$-dimensional (right) $f_4(\x)$ \eqref{equ:gaussian} with respect to the number of constraints $C$ when polynomial orders $n=2$.}
\end{figure}

\begin{figure}[htbp]
  \centering
  \label{fig:hdk}
  \includegraphics[width=.49\textwidth]{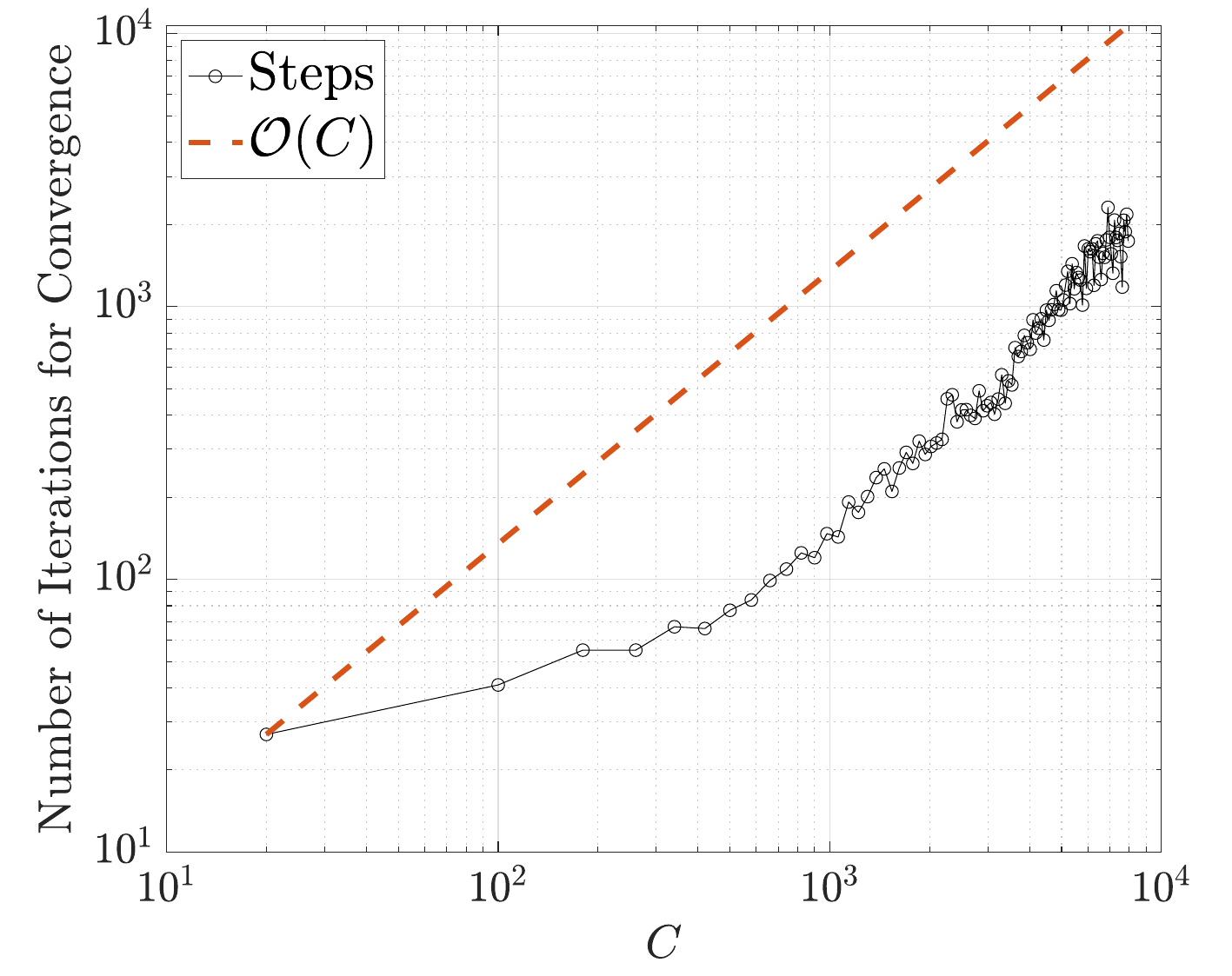}
  \includegraphics[width=.49\textwidth]{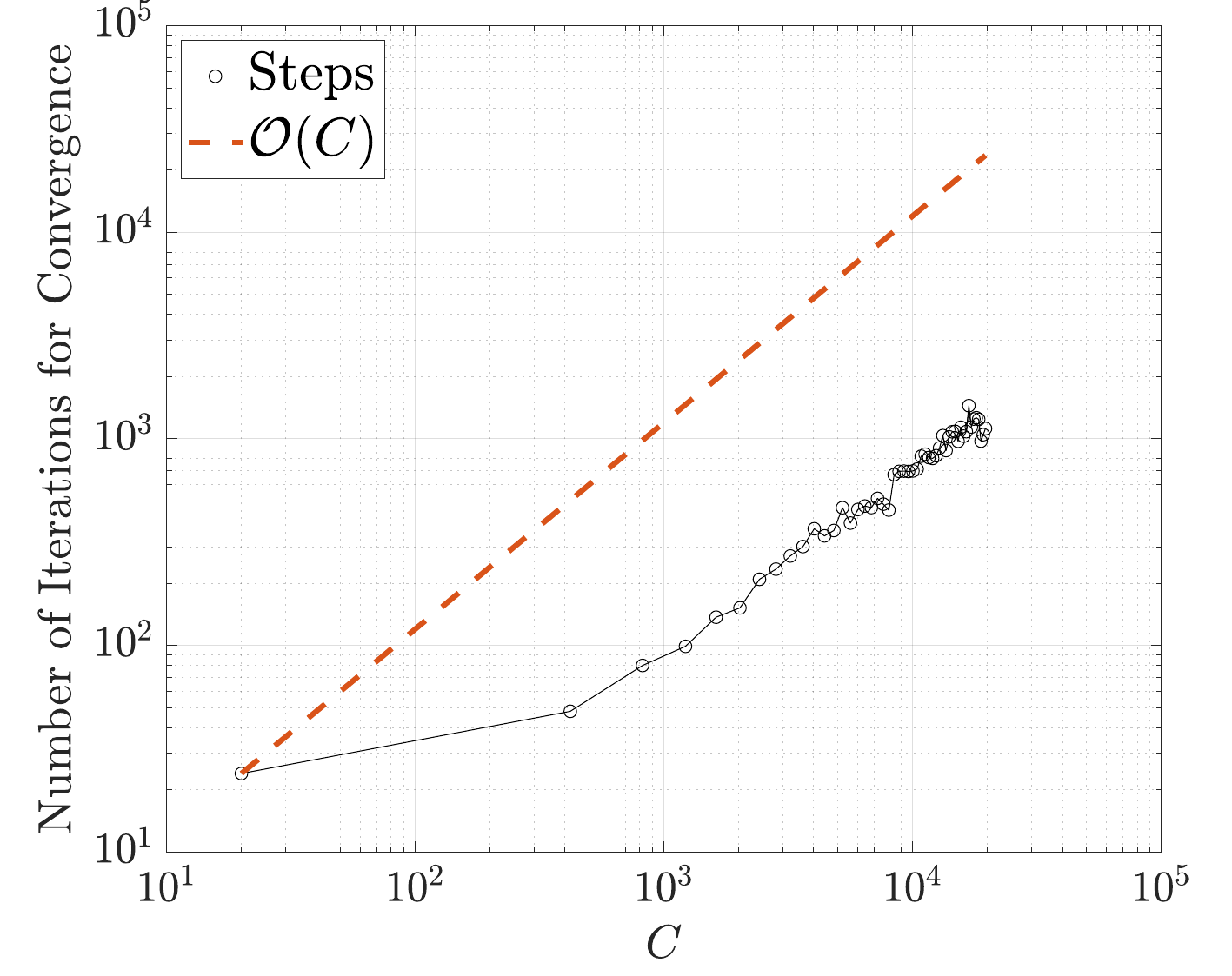}
  \caption{The number of iterations needed for restarted FISTA to converge for $100$-dimensional (left) and $200$-dimensional (right) $f_4(\x)$ \eqref{equ:gaussian} with respect to the number of constraints $C$ when polynomial orders $n=2$.}
\end{figure}

In this section, we provide the approximation results for $f_6(\x)$ at polynomial order $n=2$ in dimensions $d=100$ and $d=200$, where the cardinality of the polynomial spaces is $N=5,151$ and $N=20,301$, respectively. We use $K=3,000$ random sample points for the approximation, another $5,000$ independent random samples to evaluate the results.
The results are in Figures \ref{fig:hdapp}, \ref{fig:hdpos}, and \ref{fig:hdk}. It is observed that only a few hundreds iteration steps are needed for convergence, for the number of positivity constraint points $C\sim 10^3$ in these high-dimensional examples.

\subsubsection{Effect of Dimensionality}

\begin{figure}[htbp]
  \centering
  \label{fig:hddim}
  \includegraphics[width=.99\textwidth]{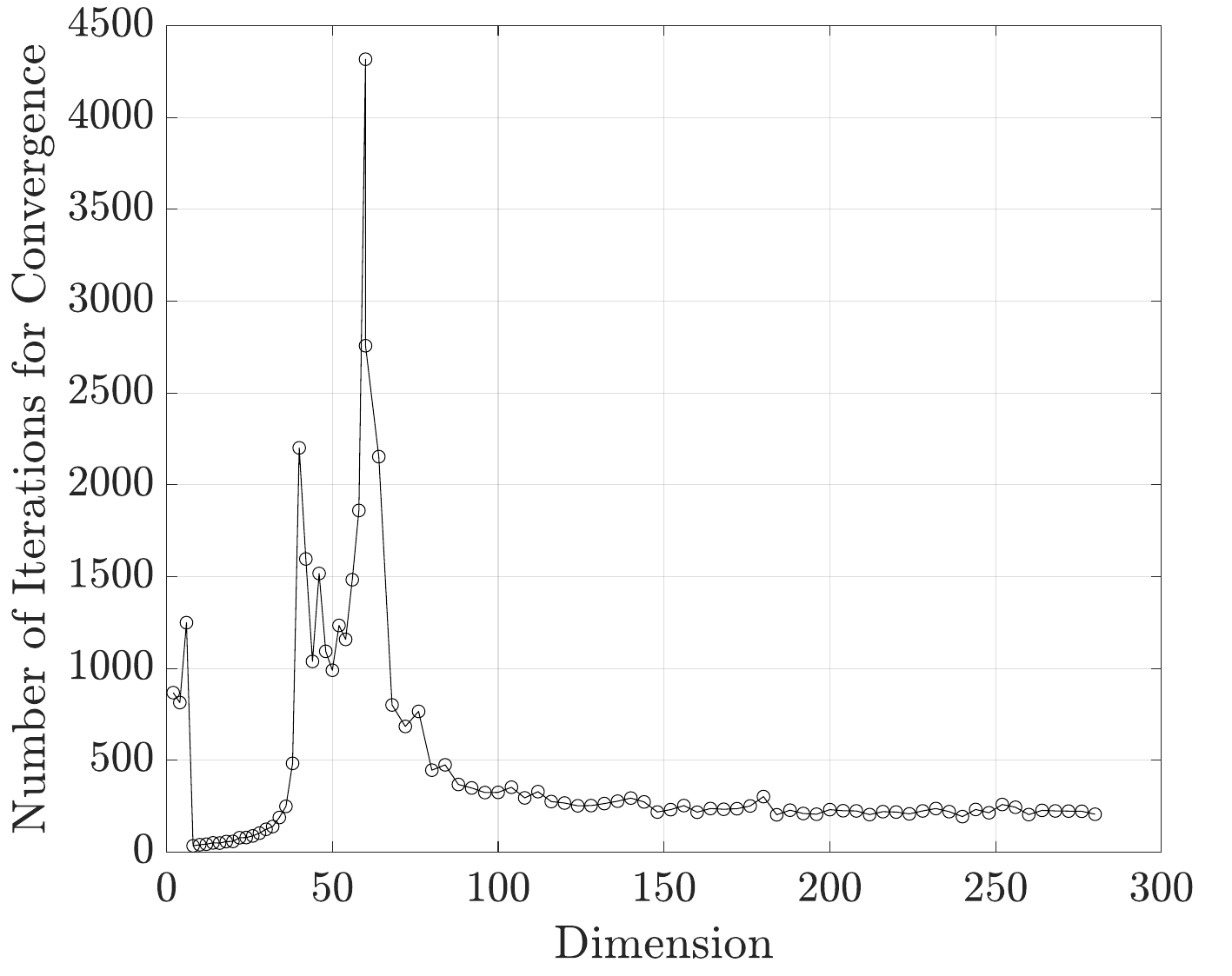}
  \caption{The number of iterations needed for restarted FISTA to converge for $f_4(\x)$ \eqref{equ:gaussian} with respect to the number of dimensions using polynomial orders $n=2$, number of approximation points $N=2,000$ and number of constraints $C=2,000$.}
\end{figure}

Finally, we examine the impact of dimensionality on the convergence. We conduct tests for $f_6(\x)$ with $2,000$ random approximation points and $2,000$ random constraint points, for dimensions from $2$ to $280$. The number of iterations needed for restarted FISTA to converge are shown in Figure \ref{fig:hddim}. \RV{We observe that, at least for this test, the number of iterations needed for the proposed FISTA with restart of high dimensions ($d \geq 100$) is significantly smaller than that of lower dimensions. This is explained by the fact that the condition of matrix $\mathbf{B}\K^\dagger\mathbf{B}^T$ is better when $d$ gets larger, with $C$ remaining unchanged. Moreover, the iteration number is not sensitive to the dimensionality beyond $d=100$.}
\section{Conclusion}
\label{sec:conclu}
In this paper, we have proposed a convex optimization based computational framework for approximating a function in high dimensions by polynomials  with non-negative and bound-preserving constraints. In particular, the restarted FISTA method applied on a proper dual problem can be easily implemented, and scales well with the problem size.  
Numerical tests have verified the effectiveness of the method for problems in a few hundred dimensions. Future work consists of exploring how to design optimal parameters scuh as the step size for Douglas-Rachford splitting to achieve faster convergence and lower cost than the restarted FISTA method. 
 
\bibliographystyle{siamplain}
\bibliography{references}

\end{document}